\documentclass[11pt]{article}
\usepackage{graphicx}
\usepackage{color}
\usepackage{amsmath}
\usepackage{amssymb}
\usepackage{amscd}
\usepackage{bbm}

\usepackage{fancyhdr,a4wide}
\usepackage{amsthm}

\newcommand{\R}{\mathbb{R}}
\newcommand{\inr}[1]{\left\langle #1 \right\rangle}

\newcommand{\E}{\mathbb{E}}
\newcommand{\M}{\mathbb{M}}
\newcommand{\Q}{\mathbb{Q}}

\newcommand{\eps}{\varepsilon}

\newcommand{\cF}{{\cal F}}

\newtheorem{Theorem}{Theorem}[section]
\newtheorem{Lemma}[Theorem]{Lemma}

\newtheorem{Definition}[Theorem]{Definition}

\newtheorem{Proposition}[Theorem]{Proposition}
\newtheorem{Corollary}[Theorem]{Corollary}

\newtheorem{Assumption}{Assumption}[section]
\newtheorem{Question}[Theorem]{Question}

\numberwithin{equation}{section}

\def \proof {\noindent {\bf Proof.}\ \ }
\def \remark {\noindent {\bf Remark.}\ \ }

\def \endproof
{{\mbox{}\nolinebreak\hfill\rule{2mm}{2mm}\par\medbreak}}
\def\IND{\mathbbm{1}}

\newcommand{\ol}{\overline}
\newcommand{\wt}{\widetilde}
\newcommand{\wh}{\widehat}
\newcommand{\argmin}{\mathop{\mathrm{argmin}}}
\newcommand{\X}{\mathcal{X}}
\newcommand{\C}{\mathcal{C}}
\newcommand{\D}{\mathcal{D}}

\newcommand{\EXP}{\mathbb{E}}
\newcommand{\PROB}{\mathbb{P}}

\newcommand{\var}{\mathrm{Var}}
\newcommand{\F}{{\mathcal F}}

\begin{document}

\title{Risk minimization by median-of-means tournaments
\thanks{
G\'abor Lugosi was supported by
the Spanish Ministry of Economy and Competitiveness,
Grant MTM2015-67304-P and FEDER, EU. Shahar Mendelson was supported in part by the Israel Science Foundation.
}
}
\author{
G\'abor Lugosi\thanks{Department of Economics and Business, Pompeu
  Fabra University, Barcelona, Spain, gabor.lugosi@upf.edu}
\thanks{ICREA, Pg. Lluís Companys 23, 08010 Barcelona, Spain}
\and
Shahar Mendelson \thanks{Department of Mathematics, Technion, I.I.T, and Mathematical Sciences Institute, The Australian National University, shahar@tx.technion.ac.il}}

\maketitle

\begin{abstract}
We consider the classical statistical learning/regression problem, when
the value of a real random variable $Y$ is to be predicted based on the
observation of another random variable $X$. Given a class of functions
$\F$ and a sample of independent copies of $(X,Y)$, one needs to
choose a function $\wh{f}$ from $\F$ such that $\wh{f}(X)$
approximates $Y$ as well as possible, in the mean-squared sense.
We introduce a new procedure, the so-called median-of-means
tournament, that achieves the optimal tradeoff between accuracy and confidence under minimal assumptions, and in particular outperforms classical methods based on empirical risk minimization.
\end{abstract}

\section{Introduction}
Estimation and prediction problems are of central importance in statistics and learning theory.
In the standard regression setup, $(X,Y)$ is a pair of random variables: $X$ takes its values in some (measurable) set $\X$ and is distributed according to an unknown probability measure $\mu$, while $Y$ is real valued that is also unknown.
Given a class $\F$ of real-valued functions defined on $\X$,
one wishes to find $f\in \F$ for which $f(X)$ is a good prediction of $Y$. Although one may consider various notions of `a good prediction', we restrict our attention to the--perhaps most commonly used--\emph{squared error}: the learner is penalized by $(f(X)-Y)^2$ for predicting $f(X)$ instead of $Y$. Thus, one would like to find a function $f \in \F$ for which the expected loss $\EXP (f(X)-Y)^2$, known as the \emph{risk}, is as small as possible. Naturally, the best performance one may hope for is of the risk minimizer in the class, that is, that of
\[
   f^*= \argmin_{f\in \F} \EXP (f(X)-Y)^2~.
\]
We  assume in what follows that the minimum is attained and $f^*\in \F$ exists and is unique, as is the case when $\F \subset L_2(\mu)$ is a closed, convex set.

One may formulate two natural goals in estimation and prediction problems. One of them is to find a function $f \in \F$ whose $L_2(\mu)$ distance to $f^*$
\begin{equation} \label{eq:estimation}
   \left(\EXP \left(f(X)-f^*(X)\right)^2 \right)^{1/2}
\end{equation}
is as small as possible. The other is to ensure that the \emph{excess risk} of the selected function
\begin{equation} \label{eq:prediction}
   R(f) = \EXP (f(X)-Y)^2 - \EXP (f^*(X)-Y)^2~,
\end{equation}
is a small.

The crucial difference between this type of problems and standard questions in approximation theory is that the available information is
limited to a random sample. One observes $\D_N=((X_1,Y_1),\ldots,(X_N,Y_N))$, that is, $N$
independent pairs, where each $(X_i,Y_i)$ has the same distribution as
$(X,Y)$ and $\D_N$ is independent of $(X,Y)$.  The fact that the
distribution of the pair $(X,Y)$ is
not known makes it impossible to invoke approximation-theoretical methods and identify directly the true minimizer of the risk.

Given a sample size $N$, a \emph{learning procedure} is a map
$\Phi:(\X \times \R)^N \to \F$ that assigns to each sample
$\D_N=(X_i,Y_i)_{i=1}^N$ a (random) function in $\F$, which we denote by
$\wh{f}_N$. The \emph{mean squared error} (also called the \emph{estimation error}) of $\Phi$ is the squared $L_2(\mu)$
distance between the true minimizer $f^*$ and the function selected by
$\Phi$ given the data $\D_N$, that is, the conditional expectation
$$
\|\wh{f}_N-f^*\|_{L_2}^2= \E \bigl((\wh{f}_N(X)-f^*(X))^2 | \D_N \bigr) \equiv {\cal E}_e^2~,
$$
where throughout the article, for $q\ge 1$, we use the notation
\[
 \|f-g\|_{L_q}=    \left(\EXP \left|f(X)-g(X)\right|^q \right)^{1/q}
\quad \text{and also} \quad
 \|f-Y\|_{L_q}=    \left(\EXP \left|f(X)-Y\right|^q \right)^{1/q}~.
\]

The \emph{excess risk}, also known as the \emph{prediction error}, compares the `predictive capabilities' of $\wh{f}_N$ to that of the best in the class, and is defined by the conditional expectation
$$
R(\wh{f}_N)=\E \bigl((\wh{f}_N(X)-Y)^2|\D_N\bigr) - \E (f^*(X)-Y)^2 \equiv {\cal E}_p^2~.
$$

Both the mean squared error and the excess risk are functions of the
given data $\D_N$, and as such are random quantities.
It is worth noting here that in the special situation when
$f^*(X)=\E(Y|X)$, we have $R(\wh{f}_N)=\|\wh{f}_N-f^*\|_{L_2}^2$.
This is the case, for example, when $Y=f(X)+W$ for some $f\in \F$
and a zero-mean random variable $W$ that is independent of $X$.
However, in general, a small mean squared error does not automatically
imply a small excess risk, or vice-versa.

In what follows we
refer to both ${\cal E}_e$ and ${\cal E}_p$ as the \emph{accuracy\footnote{Sometimes the accuracy is defined by ${\cal E}_e^2$ and ${\cal E}_p^2$.}} of learning procedure $\Phi$. The
\emph{confidence} of $\Phi$ for an error rate of ${\cal E}$ is the
probability (with respect to the product measure on $(\X \times \R)^N$
endowed by the pair $(X,Y)$) with which $\Phi$ performs with accuracy
smaller than ${\cal E}$.

Note that, up to this point, $Y$ was an arbitrary square-integrable real-valued random
variable, and obviously one would like to be able to treat as wide a
variety of targets as possible. Clearly, the accuracy and
confidence one may establish may depend on some features of the
target--for example, some a-priori estimate on its $L_q$ norm--, or on its
``distance" to $\F$, etc. We consider a broad set of
\emph{admissible targets} ${\cal Y}$, and the accuracy and confidence
of $\Phi$ relates to its performance for \emph{any} admissible target
$Y \in {\cal Y}$. Thus, a \emph{learning problem} is the triplet $(\F,Y,X)$, when $X$ and $Y$ are not known, though the learner does know that $Y \in {\cal Y}$.

It is clear that there is a tradeoff between the accuracy and
confidence in a given learning problem: the smaller the error is, the
harder it is to attain it. The question of this accuracy/confidence
tradeoff is of utmost importance in statistical learning theory, and
has been investigated extensively in numerous manuscripts since the
early days of the area in the late 1960's (see, for example, the books
\cite{VaCh74a,DGL96,vaWe96,AnBa99,vdG00,Mas06,Kolt08,TS09,BuvdG11} for a sample of the work devoted to this question). To find
the right accuracy/confidence tradeoff one must first identify a lower
bound on the tradeoff in terms of the sample size, the structure of
$\F$ and possibly some additional information on $X$ and $Y$, and then
come up with a learning procedure that attains the tradeoff.

Roughly put, one should explore the tradeoff for the set of
``achievable accuracies" of each learning problem. An accuracy ${\cal
  E}$ is achievable if there is a learning procedure in $\F$ that
achieves the accuracy ${\cal E}$ for the problem $(\F,Y,X)$ with
constant probability--say at least $3/4$--, and because $Y$ is not known, this
has to hold for any $Y \in {\cal Y}$. We define the \emph{accuracy
  edge} as the smallest achievable accuracy of a problem.

The primary question is to find the correct accuracy/confidence
tradeoff for any (reasonable) learning problem and identify a learning
procedure that attains that tradeoff all the way down to the accuracy
edge. It should be noted that up to this point in time, and other than
in a few isolated examples, no learning procedure even came close to
the optimal tradeoff at \emph{any} nontrivial accuracy level.

In this article we solve this problem by presenting an \emph{optimal
  learning procedure}: it yields the best possible accuracy/confidence
tradeoff (almost) up to the achievable edge, and under minimal
assumptions on the learning problem.

The minor reservation ``almost" is due to fact that more often than
not, the identity of the accuracy edge of a learning problem is not
known. As it is explained in what follows, while one may provide
lower estimates on the accuracy edge, there is a very real possibility
that such estimates are too optimistic and the real accuracy edge is
larger. Regardless, the procedure we introduce ``gets as close" to the
accuracy edge as any other known procedure--or better, and with a
dramatically better confidence. It also exhibits the optimal
accuracy/confidence tradeoff for larger errors, something that no
other procedure is known to do.

It should be clarified at this point that
by ``optimal accuracy'' we mean an accuracy that is optimal up to a constant factor
and optimality in the confidence means that the probability with which
the claimed accuracy does not hold is optimal up to a constant factor
in the exponent.

Before we present a more accurate formulation of our main results and
describe the optimal procedure, let us explain what our procedure is
not: it is not \emph{empirical risk minimization} (ERM), nor any of
its ``family members".

\subsection{ERM--the wrong choice}
Perhaps the most natural way of choosing $\wh{f}_N$ is by empirical risk minimization, that is, by \emph{least squares regression},
\[
  \wh{f}_N = \argmin_{f\in \F} \sum_{i=1}^N (f(X_i)-Y_i)^2~.
\]
Again, we assume that the minimum is attained, while if there are several minimizers, $\wh{f}_N$ may be chosen among them in an arbitrary way.

The performance of least squares regression has been thoroughly studied in many different scenarios. A sample of the rich literature includes Gy\"orfi, Kohler, Krzyzak, Walk \cite{GyKoKrWa02}, van de Geer \cite{vdG00}, Bartlett, Bousquet, and Mendelson \cite{BaBoMe05}, Koltchinskii \cite{Kol11}, Massart \cite{Mas06}.

The simple idea behind empirical risk minimization is that, for each
$f\in \F$, the empirical risk $(1/N) \sum_{i=1}^N (f(X_i)-Y_i)^2$ is a
good estimate of the risk $\EXP (f(X)-Y)^2$ and the minimizer of the
empirical risk should nearly match that of the ``true'' risk.
%Indeed,
%the performance of the empirical risk minimizer depends on the
%fluctuations of the empirical risk around its expected value near the
%minimizer, and the latter depends on the geometry of the class $\F$.

It turns out that the performance of ERM changes dramatically
according to the tail behaviour of the functions involved in the given
learning problem. One may show (see, e.g., \cite{LeMe16}) that if $\F$
is convex and the functions in $\F$ (more precisely, the random
variables $\{f(X) : f \in \F\}$) and the target $Y$ have well-behaved
tails (and by ``well-behaved" we mean sub-Gaussian), ERM preformed in
$\F$ yields good results: for an accuracy that is not far from the
accuracy edge, it attains the optimal confidence, though it does not
maintain the optimal accuracy/confidence tradeoff for larger
errors. Unfortunately, the situation deteriorates considerably when
either members of $\F$ or one of the admissible targets is
heavy-tailed in some sense. In such cases, the performance of ERM is
significantly weaker than the known theoretical limitations of the
accuracy/confidence tradeoff. Moreover, replacing ERM with a
different procedure is of little use: other than in few and rather
special learning problems, there have been no known alternatives to
ERM whose performance comes close to the known theoretical limitations
of the accuracy/confidence tradeoff, and certainly not when the
problem is heavy-tailed.

The reason for ERM's
diminished capacity is that it is sensitive to even a small number of
atypical points in the sample $(X_i,Y_i)_{i=1}^N$: since ERM
selects a minimizer of the empirical mean of the squared loss,
atypical values may distort the selection and send ERM to the wrong
part of $\F$. This sensitivity is clearly reflected in the confidence
with which ERM operates in heavy-tailed situations: roughly put, one
can guarantee that ERM performs with the right accuracy only on
samples that are not contaminated by a significant number of atypical
values. However, in heavy-tailed situations, the latter does not occur
frequently, and having atypical values is simply a fact of life one has to deal with.

In contrast, the procedure we suggest as an alternative to ERM leads to the optimal accuracy/confidence tradeoff even in heavy-tailed situations. Unlike ERM, it is not sensitive even to a large number of atypical sample points.

Before we dive into a more technical description of our results, let us present the following classical example of linear regression in $\R^n$, exhibiting
the limitations of empirical risk minimization in heavy-tailed problems, and comparing its performance to that of the procedure we introduce.

Let $\F=\{\inr{t,\cdot} : t \in \R^n\}$ be the class of linear functionals on $\R^n$. Let $X$ be an isotropic random vector in $\R^n$ (i.e., $\E\inr{t,X}^2 = 1$ for every $t$ in the Euclidean unit sphere) and assume that $X$ exhibits some (very weak) norm equivalence in the following sense: there are $q>2$ and $L>1$ for which, for every $t \in \R^n$, $\|\inr{X,t}\|_{L_q} \leq L \|\inr{X,t}\|_{L_2}$.

Assume that one is given $N$ noisy measurements of $\inr{t_0,\cdot}$
for a fixed but unknown $t_0 \in \R^n$. Specifically, assume that
$Y=\inr{t_0,X}+W$ for some
symmetric random variable $W$ that is independent of $X$ and has
variance $\sigma^2$.
One observes the ``noisy" data $(X_i,Y_i)_{i=1}^N$ and
the aim is to approximate $t_0$
with a small error (accuracy) and with high probability (confidence)
using this random data only.

One may show
that a nontrivial estimate is possible only when $N \geq
cn$ for a suitable absolute constant $c$, and we consider only such
values of $N$. Also, there are known estimates on the theoretical limitations of
this problem: a lower bound on the accuracy edge is of the order of $\sigma
\sqrt{n/N}$, and for an accuracy level that is proportional to the accuracy edge, say, $c_0\sigma\sqrt{n/N}$ for a suitable absolute constant $c_0$, the conjectured confidence is $1-2\exp(-c_1n)$.

If there is no information on
higher than $q$-th moments for linear functionals, and no
information beyond the second moment for $W$ is available,
this clearly is a (potentially) heavy-tailed scenario.
It turns out (the claims made here follow from results of \cite{LeMe16},
see the next section for the general statements)
that the
best that one can guarantee using ERM is a choice of $\wh{t} \in
\R^n$, for which the Euclidean norm $\|\wh{t}-t_0\|_2 = \|\inr{\wh{t},X} -
\inr{t_0,X}\|_{L_2} \leq r$ with probability at least
$1-\delta-2\exp(-cN)$; the error $r$ is defined as the smallest number
for which
$$
(*)=\left\|\frac{1}{\sqrt{N}}\sum_{i=1}^N W_i X_i\right\|_2 \leq c\sqrt{N} r \ \ \ {\rm with \ probability \ at \ least \ } 1-\delta~.
$$
Since $X$ is isotropic, one has $\E\|X\|_2^2 = n$.
Therefore, the mean of (*) is bounded as
$$
\E\left\|\frac{1}{\sqrt{N}}\sum_{i=1}^N W_i X_i\right\|_2
\leq \left(\E\left\|\frac{1}{\sqrt{N}}\sum_{i=1}^N W_i X_i\right\|_2^2\right)^{1/2} = \sigma (\E\|X\|_2^2)^{1/2} =\sigma \sqrt{n}~.
$$
Moreover, because of the minimal assumptions on $W$ and $X$, the best estimate one can hope for on $(*)$ and that holds with probability $1-\delta$ follows from Chebyshev's inequality. In particular, it leads to the following, rather unsatisfactory, estimate on the performance of ERM:
$$
\|\wh{t}-t_0\|_2 \leq \frac{c_0(q,L) }{\delta} \sigma\sqrt{\frac{n}{N}} \ \ \ {\rm with \ probability \ } 1-\delta-2\exp(-c_1 N)~.
$$
Also, if one wishes for an error that is proportional to the (conjectured) accuracy edge, that is, of the order of $\sigma \sqrt{n/N}$, the best that one can hope for is a constant confidence--a very different estimate from the conjectured confidence of $1-2\exp(-c_1n)$.

Although what we have described is an upper bound, one may show (see
\cite{LeMe16}) that this estimate captures the performance of empirical risk minimization, and in particular exhibits ERM's inability to deal with atypical
sample points. The reality is that ERM performs with an accuracy of
the order of $\sigma \sqrt{n/N}$ only on the relatively few samples
that contain almost no misleading data.

%On the other hand, known estimates on the theoretical limitations of
%this problem are quite different. A lower bound on the accuracy edge is $c\sigma
%\sqrt{n/N}$, and the conjectured confidence at that level is $1-2\exp(-c_1n)$, which is
%significantly better than the estimate of constant probability when
%ERM is used to select $\wh{t}$.

The main result of this article, when applied to this example, shows that under the same assumptions, the procedure we suggest selects $\wh{t}$ for which
\begin{equation} \label{eq:reg-r-n}
\|\wh{t}-t_0\|_2 \leq C \sigma\sqrt{\frac{n}{N}} \ \ \ {\rm with \ probability \ } 1-2\exp(-cn)
\end{equation}
for some numerical constants $c,C>0$; that is, it performs with optimal confidence at a level that is proportional to the accuracy edge. In fact, our procedure gives the optimal confidence for any accuracy $r \geq c^\prime\sigma \sqrt{n/N}$.

Note that for the special case of linear regression described above, Hsu and Sabato
\cite{HsSa16} achieve slightly (by a factor logarithmic in $n$) weaker bounds than \eqref{eq:reg-r-n} under slightly stronger ($(4+\epsilon)$-th moment) assumptions. We also refer to Minsker
\cite{Min15} for related bounds for sparse regression under possibly heavy-tailed
variables. Of course, these results hold in a rather special example, while our main result yields optimal estimates for almost any convex class $\F$ and target $Y$, and not just for linear regression in $\R^n$.

In the next section we present the required definitions, outline the
current state of the art, and formulate our main results.
In Section \ref{sec:tournament}, we  describe the new procedure in
detail.
In Section \ref{sec:examples} we illustrate the power of the main
results on some canonical examples,
before turning to the proofs of our results.

Let us point out the well-understood fact that the behaviour of the
accuracy and confidence in learning problems in which $\F$ is not
convex is trivial in some sense, and totally different from the convex
case, which is why focus on the latter. As it happens, the dominating
factor in non-convex problems is the `location' of the targets $Y$
relative to $\F$ rather than the structure of $\F$, and an
`unfavourable location' of a target $Y$ completely distorts the
accuracy and confidence of the learning problem. However, even
if $\F$ is not convex, all the targets of the form $Y=f_0(X)+W$ for
$f_0 \in \F$ and $W$ that is symmetric and independent of $X$ happen
to be in a `favourable' location and thus our results apply to such
problems as well.

%HERE WE NEED REFERENCES TO MINSKER, HSU and SABATO, AUDIBERT and CATONI

\section{The accuracy edge and the accuracy/confidence tradeoff}

We begin by describing the known theoretical limitations on the accuracy edge and on the
accuracy/confidence tradeoff for a given learning problem. To this end, let us introduce some notation, following the path of \cite{Men15,Men16}.

Let $D=\{f: \|f\|_{L_2}\le 1\}$ be the unit ball in $L_2(\mu)$ and set $S
=\{f: \|f\|_{L_2}= 1\}$ to be the unit sphere.  For
$h \in L_2(\mu)$ and $r>0$, put $rD_h = \{f : \|f-h\|_{L_2} \leq r\}$. Let
$$
{\rm star}(\F,h) = \{ \lambda f + (1-\lambda)h \ : \ 0 \leq \lambda
\leq 1, \ f \in \F\}~.
$$
Thus, ${\rm star}(\F,h)$ is the star-shaped hull of $\F$ around
$h$, that is, the union of all segments for which one end-point
is $h$ and the other is in $\F$.

The star-shaped hull ${\rm star}(\F,h)$ adds regularity to the class
around the fixed centre $h$: on the one hand, it does not increase the
size of the class by much, while on the other hand, it implies that every
function of the form $f-h$ has a `scaled-down' version when one moves towards $0$. In particular, the level sets ${\rm star}(\cF-h,0) \cap rS$ become `richer' as $r$ gets smaller: each one of them contains scaled-down copies of all `higher' levels.

Consider the localization of $\F$
$$
\F_{h,r} = {\rm star}(\F-h,0) \cap
rD,
$$
which is given by the shift that maps the designated point $h$ to
$0$. Then the resulting class is made `more regular' by taking its
star-shaped hull around $0$, and finally it is localized, by
considering its intersection with $rD$, the $L_2(\mu)$ ball of radius
$r$, centred in $0$.

Observe that if $\F$ is convex then
for any $h \in \F$, ${\rm star}(\F,h)=\F$. Also, in that case,
$$
\F_{h,r}=\{f-h : f \in \F, \ \|f-h\|_{L_2} \leq r\}~,
$$
and if, in addition, $\F$ is centrally symmetric (that is, if $f \in \F$ then $-f \in \F$), then $\F_{h,r} \subset 2\F \cap rD$.

One way of deriving lower estimates on the accuracy edge and on the accuracy/confidence tradeoff is based on the packing numbers of the localizations $\F_{h,r}$.

\begin{Definition} \label{def:pack-and-cover}
Given a set $H \subset L_2(\mu)$ and $\eps>0$, denote by
${\cal M}(H,\eps D)$ the cardinality of
a maximal $\eps$-separated subset of $H$.  That is, ${\cal M}(H,\eps
D)$  is the maximal cardinality of a subset $\{h_1,...,h_m\} \subset H$, for which $\|h_i-h_j\|_{L_2} \geq \eps$ for every $i \not = j$.
\end{Definition}

Note that if $H^\prime$ is a maximal $\eps$-separated subset of $H$,
then it is also an $\eps$-cover of $H$ in the sense that for every $h \in H$ there is
some $h^\prime \in H^\prime$ that satisfies $\|h^\prime -h\|_{L_2} \leq \eps$.

\begin{Definition} \label{def-lambda-Q}
For $\kappa,\eta>0$ and $h \in \F$, set
\begin{equation} \label{eq:covering-fixed-point-1}
\lambda_{\Q}(\kappa,\eta,h) =  \inf\{ r: \log {\cal M}(\F_{h,r},\eta r D) \leq \kappa^2 N\}~,
\end{equation}
and let
\[
\lambda_{\Q}(\kappa,\eta) = \sup_{h \in \F}
\lambda_{\Q}(\kappa,\eta,h)~.
\]
\end{Definition}

For every fixed $h \in F$, the parameter $\lambda_{\Q}$ pin-points the
level $r$ at which the localization $\F_{h,r}$ becomes ``too rich" in
the following sense: given the sample size $N$, $\F_{h,r}$ contains a
subset of cardinality $\exp(\kappa^2 N)$ that is $\eta r$-separated
with respect to the $L_2(\mu)$ norm. Note that $\lambda_{\Q}$ is not
affected by the fine structure of $\F_{h,r}$. Indeed, the set
$\F_{h,r} \cap (\eta r/2)D$ cannot contain more than two points that are $\eta r$-separated, and thus it does not contribute to the existence of a large $\eta r$-separated set in $\F_{h,r}$.

One may show that $\lambda_{\Q}$ serves as a lower bound on the accuracy edge of a learning problem, when $\F$ is convex and centrally symmetric and the admissible targets are noise-free: that is, ${\cal Y}=\{f(X) : f \in \F\}$.

\begin{Proposition} \label{thm:lambda-Q-lower} \cite{Men16}
There exist absolute constants $\kappa$ and $\eta$ for which the
following holds. Let $\F \subset L_2(\mu)$ be convex and centrally symmetric. For any learning procedure $\Phi$ there exists an $f_0\in \cF$ and target $Y=f_0(X)$ for which, with probability at least $1/4$,
$$
  \|\Phi(\D_N)-f_0\|_{L_2}^2 \geq
  \lambda_{\Q}(\kappa,\eta)~.
$$
\end{Proposition}

The following variant of $\lambda_{\Q}$ also serves as a lower bound on the accuracy edge, this time, because of `noisy' targets.

\begin{Definition} \label{def-lambda-M}
For $\kappa>0$, $0<\eta<1$ and $h \in \F$, set
\begin{equation} \label{eq:covering-fixed-point-2}
\lambda_{\M}(\kappa,\eta,h) =  \inf\{ r: \log {\cal M}(\F_{h,r},\eta r D) \leq \kappa^2 N r^2\}
\end{equation}
and let
\[
\lambda_{\M}(\kappa,\eta) = \sup_{h \in \F}
\lambda_{\M}(\kappa,\eta,h)~.
\]
\end{Definition}

\begin{Proposition} \label{thm:lambda-M-lower} (See, e.g., \cite{Men16}.)
There exist absolute constants $\kappa$ and $\eta$ for which the
following holds.
Let $\F \subset L_2(\mu)$ and set $W$ to be a centred Gaussian
variable with variance $\sigma>0$ that is independent of $X$.
Then, for any learning procedure $\Phi$ there exists $f_0\in \F$ and a target $Y=f_0(X)+W$  for which, with probability at least $1/4$,
$$
\|\Phi(\D_N)-f_0\|_{L_2}^2 \geq \lambda_{\M}(\kappa/\sigma,\eta)~.
$$
\end{Proposition}
In particular, if ${\cal Y}$ contains all the targets of the form $f_0(X)+W$, for $f_0 \in \F$ and $W$ that is a centred Gaussian variable with variance $\sigma$ that is independent of $X$, then the accuracy edge is at least $\lambda_{\M}(\kappa/\sigma,\eta)$.

Combining these two facts, we have a lower bound on the accuracy edge:
\[
%begin{equation} \label{eq:accuracy-edge}
\lambda^* \equiv \max\{\lambda_{\Q}(\kappa_1,\eta_1),\lambda_{\M}(\kappa_2/\sigma,\eta_2)\}~,
\]
%end{equation}
for some constants $\kappa_i,\eta_i$, $i=1,2$. However, there is no guarantee that this lower estimate is sharp. As we explain in what follows, it is very possible that the true accuracy edge is larger.

Let us turn to the theoretical limitations of the accuracy/confidence tradeoff, assuming that the set of admissible targets ${\cal Y}$ is not too trivial. By that we mean that it at least
contains all the targets of the form $Y=f_0(X)$ (the noise-free problems) and $Y=f_0(X)+W$, for $f_0 \in \F$ and $W$ that is a centred Gaussian variable with variance $\sigma$,
and is independent of $X$. We call this set of targets \emph{minimal}, and of course, ${\cal Y}$ could be much larger.

Applying the results in \cite{LeMe16} and \cite{Men16} one has the following:
\begin{Proposition} \label{thm:lower-general-r}
There exists an absolute constant $c$ for which the following holds. Let $\F \subset L_2(\mu)$ be a class that is star-shaped around one of its points (i.e., for some $f_0\in\F$ and every $f\in\F$, $[f_0,f] \subset \F$). Consider ${\cal Y}$ that contains the minimal set of targets. If $\Phi$ is a learning procedure that performs with accuracy $r$ and confidence $1-\delta$ for every such target, then
$$
\delta \geq \exp(-cN\min\{1,r^2/\sigma^2\}).
$$
\end{Proposition}

These facts set our first benchmark (which may be too optimistic, of course): the lower bound on the accuracy edge
\begin{equation} \label{eq:optimal-accuracy}
\lambda^* =\max\{\lambda_{\Q}(\kappa_1,\eta_1),\lambda_{\M}(\kappa_2/\sigma,\eta_2)\}~,
\end{equation}
and the bound on the accuracy/confidence tradeoff for $r \geq c_0 \lambda^*$,
\begin{equation} \label{eq:optimal-confidence}
1-2\exp \bigl(-c_1 N\min\{1,\sigma^{-2}r^2\}\bigr)~.
\end{equation}

As we noted earlier, $\lambda^*$ is an optimistic, and perhaps not very
realistic, lower bound on the accuracy edge. A more reasonable
conjecture relies on more ``global" parameters that take into account
the fine structure of $\F$ at an arbitrarily small level, defined next.

From here on, let $(\eps_i)_{i=1}^N$ be independent, symmetric $\{-1,1\}$-valued random variables that are independent of $(X_i,Y_i)_{i=1}^N$.
\begin{Definition} \label{def:fixed-emp}
For $\kappa>0$ and $h \in \F$ let
\begin{equation} \label{eq:emp-fixed-point}
r_{E}(\kappa,h) = \inf\left\{ r: \E \sup_{u \in \F_{h,r}} \left| \frac{1}{\sqrt{N}} \sum_{i=1}^N \eps_i u(X_i) \right| \leq \kappa \sqrt{N} r\right\},
\end{equation}
and set $r_{E}(\kappa) = \sup_{h \in \F} r_{E}(\kappa,h)$.
\end{Definition}

The parameter $r_E(\kappa,h)$ measures the empirical oscillation
around $h$. It does not depend on the identity of the target $Y$ and
is a purely intrinsic parameter of the class $\F$. However, it may be
highly affected by functions in $\F$ that are close to $h$, and as
such it is more ``global" than $\lambda_{\Q}$.

The other ``global" parameter we require does depend on $Y$. It is used to calibrate the interaction between $\F$ and the target.

\begin{Definition} \label{def:fixed-multi}
For $\kappa>0$ and $h \in \F$ set $\ol{r}_{\M}(\kappa,h)$ to be
\begin{equation}
\ol{r}_{\M}(\kappa,h) = \inf\left\{ r: \E \sup_{u \in \F_{h,r}} \left| \frac{1}{\sqrt{N}} \sum_{i=1}^N\eps_i u(X_i) \cdot (h(X_i)-Y_i) \right| \leq \kappa \sqrt{N} r^2\right\}.
\end{equation}
For $\sigma>0$ put $\F_Y^{(\sigma)}=\{ f \in \F : \|f(X)-Y\|_{L_2} \leq \sigma\}$ and let
${\wt r}_{\M}(\kappa,\sigma) = \sup_{h \in \F_Y^{(\sigma)}} {\ol r}_{\M}(\kappa,h)$.
\end{Definition}

\remark The role of $\sigma$ and of $\F_Y^{(\sigma)}$ in Definition
\ref{def:fixed-multi} deserves some explanation. The mean oscillation
in Definition \ref{def:fixed-multi} involves the ``multipliers''
$(h(X_i)-Y_i)_{i=1}^N$, and the right choice for
the centre $h$ is the unknown $f^*$. While one may take into account
the ``worst" $h \in \F$, doing so makes little sense, as $f^*$ is
the minimizer of the $L_2$ distance between $\F$ and $Y$, and for the
worst $h$, $\|h-Y\|_{L_2}$ could be significantly larger than $\|f^*-Y\|_{L_2}$. To
overcome this obstacle we assume that
an a-priori estimate on the $L_2$ distance between $Y$ and $\F$
(i.e.,  a value $\sigma$ such that $\|f^*-Y\|_{L_2} \leq \sigma$)
is available. With this information one still
needs to consider the worst centre $h$, but only among all functions
$h \in \F$ that satisfy $\|h-Y\|_{L_2} \leq \sigma$. As we explain in
Section \ref{sec:examples}, thanks to known estimates for the expectation of the supremum of a multiplier process, one only needs to keep in mind that the multipliers $\xi_i=h(X_i)-Y_i$ are
independent copies of some random variable $\xi$ that satisfies some
moment condition, such as
$\|\xi\|_{L_q} \leq L\sigma$ for some $q>2$ and a suitable constant $L$.
%For example, one may show (see, e.g., \cite{VW}) that if $\xi$ is independent of $X$ then
%$$
%\E \sup_{h \in H} \left|\sum_{i=1}^N \eps_i \xi_i h(X_i) \right| \leq C \|\xi\|_{L_{2,1}} \E\sup_{h \in H} \left|\sum_{i=1}^N \eps_i h(X_i) \right|~,
%$$
%where $\|\xi\|_{L_{2,1}} = \int_0^\infty \PROB^{1/2}(|\xi| > t)dt$, and the right-hand side may now be handled by standard empirical process techniques.

In light of the results from \cite{Men15}, a realistic alternative to $\lambda^*$ is
\begin{equation} \label{eq:r*}
r^*=\max\{\lambda_{\Q}(c_1,c_2),\lambda_{\M}(c_1/\sigma,c_2),r_E(c_1),\wt{r}_{\M}(c_1,\sigma) \}~,
\end{equation}
for some constants $c_1,c_2$, and when for the given (and unknown)
target $Y \in {\cal Y}$ one has $\|Y-f^*(X)\|_{L_2} \leq \sigma$. Indeed, in
\cite{Men15} it was shown that under some mild conditions on the
learning problem, specified below, ERM performs in $\F$
with accuracy $cr^*$ and constant confidence. Thus, $r^*$ is a potential (and to-date, the best)
candidate for the accuracy edge. 

Therefore, up to the issue of the true
identity of the accuracy edge, the (somewhat vaguely formulated) question of the accuracy/confidence tradeoff is as follows:
\begin{Question} \label{qu:tradeoff}
Is there a learning procedure which, for any reasonable learning problem and any $r \geq c_1r^*$, performs with accuracy $r$ and confidence $1-2\exp \bigl(-c_2 N\min\{1,\sigma^{-2}r^2\}\bigr)$, thus achieving the optimal accuracy/confidence tradeoff?
\end{Question}

Our main result answers Question \ref{qu:tradeoff} in the affirmative
where our notion of a ``reasonable learning problem'' is formulated
rigorously below.

The first theorem we present requires the following conditions:
\begin{Assumption} \label{ass:moment}
Let $L$ be a constant and let $X$ be distributed according to the measure $\mu$ on ${\cal X}$.
Given a locally compact, convex class of functions $\F \subset L_2(\mu)$ and $Y \in L_2$, assume that
\begin{description}
\item{$\bullet$} for every $f,h \in \F$, $\|f-h\|_{L_4} \leq L \|f-h\|_{L_2}$;
\item{$\bullet$} for every $f \in \F$, $\|f-Y\|_{L_4} \leq L\|f-Y\|_{L_2}$;
\item{$\bullet$} $\|f^*-Y\|_{L_2} \leq \sigma$ for some known
  constant $\sigma >0$.
\end{description}
\end{Assumption}

\begin{Theorem} \label{thm:main-general}
Let $L \geq 1$, $\sigma>0$, and suppose Assumption \ref{ass:moment}.
There exist constants $c,c_0,c_1$ and $c_2$ that
depend only on $L$ for which the following holds.
Let
$$
r^*(f^*)=\max\{\lambda_{\Q}(c_1,c_2,f^*),\lambda_{\M}(c_1/\sigma,c_2,f^*),r_E(c_1,f^*),\ol{r}_{\M}(c_1,f^*) \}~,
$$
and fix $r \geq 2r^*(f^*)$.

There exists a procedure that, based on the data
$\D_N=(X_i,Y_i)_{i=1}^N$ and the values of $L$, $\sigma$ and $r$, selects a
function $\wh{f}\in \F$ such that, with probability at least
$$
1-\exp\bigl(-c_0N \min\{1,\sigma^{-2}r^2\}\bigr)~,
$$
\[
\|\wh{f}-f^*\|_{L_2} \leq c r \ \ {\rm  and \ \ }
\E \bigl((\wh{f}(X)-Y)^2 |\D_N \bigr) \leq \E (f^*(X)-Y)^2 + (c r)^2~.
\]
\end{Theorem}

Of course, the identity of $f^*$ is not known, and therefore, it is
not reasonable to expect that $r^*(f^*)$ is known beforehand.
A ``legal'' data-independent choice
is any $r \geq 2r^*$ that is larger than $2r^*(f^*)$ regardless of the
identity of $f^*$. In
particular, Theorem \ref{thm:main-general} gives the optimal
accuracy/confidence tradeoff for any accuracy $r \geq 2r^*$.
Alternatively, one may consider the ``right'' choice of the parameter
$r$ as a model selection problem that may be selected using cross
validation if independent data are available. We do not discuss the
rather straightforward details further here.

\medskip

\remark In our main assumption of Theorem
\ref{thm:main-general}
we use the equivalence between the $L_4$ and $L_2$ norms for functions of the form $Y-f(X)$.
This allows us to derive bounds in terms of the variance of $Y-f^*(X)$.
In fact, if instead of norm equivalence we just have a bound on $\sigma_4= \|Y-f^*(X)\|_{L_4}$,
the arguments work equally well, with $\sigma_4$ replacing $\sigma$.  The sets $\F_Y^{\sigma}$ need to be adjusted as well: it should be replaced by all functions in $\F$ whose $L_4$ distance to $Y$ is at most $\sigma_4$.

\medskip

It turns out that, when dealing with independent noise, the assumptions
required in
Theorem \ref{thm:main-general} may be relaxed even further. In
particular, we do not require convexity of the class $\F$ and the
assumption of norm equivalence may be relaxed:

\begin{Theorem} \label{thm:main-independent}
Let $q>2$, $L>1$, and
  $\sigma>0$. There exist constants $c,c_0,c_1$ and $c_2$ that depend
  only on $q$ and $L$ for which the following holds. Let $\F$ be a locally compact 
  class of functions and assume that for every $f \in {\rm span}(\F)$,
  $\|f\|_{L_q} \leq L\|f\|_{L_2}$. Assume further that $Y=f_0(X)+W$
  for some $f_0 \in \F$ and $W$ that is mean-zero, independent of
  $X$, and satisfies $\|W\|_{L_2} \leq \sigma$.

Let $r^*(f^*)$ be as above and fix $r \geq 2r^*(f^*)$. There exists a procedure that, based on the data
$\D_N=(X_i,Y_i)_{i=1}^N$ and the values of $L,q$, $\sigma$ and $r$, selects a
function $\wh{f}\in \F$ such that, with probability at least
$$
1-\exp\bigl(-c_0N \min\{1,\sigma^{-2}r^2\}\bigr)~,
$$
$$
\|\wh{f}-f_0\|_{L_2} \leq c r \ \ {\rm and} \ \ \E
\bigl((\wh{f}(X)-Y)^2 |\D_N \bigr) \leq \E (f_0(X)-Y)^2 + (c r)^2~.
$$
\end{Theorem}

Note that in the case of independent additive noise, the assumptions
of Theorem \ref{thm:main-independent} are almost the minimal needed
for the learning problem to be well defined: norm equivalence for $q$
that may be arbitrarily close to $2$ and $W \in L_2$ that perhaps does
not have any higher moments. Even under these minimal assumptions, we
still obtain the optimal accuracy/confidence tradeoff.

\vskip0.4cm

\subsection{State of the art}
To put Theorem \ref{thm:main-general} in perspective, we describe the
current sharpest estimates on the accuracy and confidence of a
learning problem, focusing on possibly heavy-tailed distributions.

Firstly, there were no known results that are based on the ``averaged'' parameter $\ol{r}_{\M}$. Instead, the interaction between class members and the target were measured using the following ``in-probability" version of $\ol{r}_{\M}$:

\begin{Definition} \label{def:r-m}
For every $\kappa>0$, $0 \leq \delta \leq 1$ and $h \in \F$, set
$r_{\M}(\kappa,\delta,h)$ to be the infimum of the set of all values of $r$ for which
$$
\PROB\left(\sup_{u \in \F_{h,r}} \left| \frac{1}{\sqrt{N}} \sum_{i=1}^N\eps_i u(X_i) \cdot (h(X_i)-Y_i) \right| \leq \kappa \sqrt{N} r^2\right) \geq 1-\delta~.
$$
\end{Definition}

The best known estimate on prediction and estimation in a general convex class are based on a quite weak condition, rather than the norm equivalence we use. Recall that a class $\F$ satisfies a small-ball condition with constants $\kappa_0$ and $\rho_0$ if for every $f,h \in \F \cup \{0\}$,
$$
\PROB(|f-h| \geq \kappa_0 \|f-h\|_{L_2} ) \geq \rho_0~.
$$

\begin{Theorem} \label{thm:best-known} (Mendelson \cite{Men15}.)
Let $\F \subset L_2(\mu)$ be a convex class that satisfies the
small-ball condition with constants $\kappa_0$ and $\rho_0$, and let
$Y \in L_2$. If $r = \max\{r_E(c_1),r_{\M}(c_2,\delta,f^*)\}$ and
$\wh{f}$ is selected in $\F$ using empirical risk minimization, then, with probability at least
\begin{equation} \label{eq:old-confidence}
1-\delta-2\exp(-c_3N)~,
\end{equation}
\begin{equation} \label{eq:old-accuracy}
\|\wh{f}-f^*\|_{L_2} \leq r \ \ \ {\rm and} \ \ \ \E\bigl((\wh{f}(X)-Y)^2|\D_N\bigr) \leq \E(f^*(X)-Y)^2 + r^2,
\end{equation}
for constants $c_1,c_2$ and $c_3$ that depend only on $\kappa_0$ and $\rho_0$.
\end{Theorem}

The obvious weakness of Theorem \ref{thm:best-known} is the poor
tradeoff between the accuracy term $r_{\M}(c_2,\delta,f^*)$ and the
confidence $\delta$. There is no hope of obtaining a high confidence
result--say as in \eqref{eq:optimal-confidence}--, unless both $\F$ and $Y-f^*(X)$ exhibit a sub-Gaussian tail behaviour. If not, then for $\delta$ as in \eqref{eq:optimal-confidence}, the value  $r_{\M}(c_2,\delta,f^*)$ is very large, and the resulting accuracy estimate is rather useless--far worse than $r^*$. Moreover, replacing the small-ball assumption with some norm equivalence as in Theorem \ref{thm:main-general} does not improve the outcome. Thus, Theorem \ref{thm:best-known} is significantly weaker than Theorem \ref{thm:main-general} in every aspect.
This phenomenon exhibits the nature of ERM: it does not perform with both high accuracy \emph{and} high confidence in heavy-tailed situations, and falls well short of our benchmarks, but still, it was the best available alternative prior to this work.

\subsection*{Regression in $\R^n$ revisited}

We now show how our general results imply the optimal accuracy/confidence tradeoff for
linear regression in $\R^n$ all the way to a number proportional to the conjectured accuracy edge.

Recall that the class of functions in question is $\F=\{\inr{t,\cdot}
: t \in \R^n\}$,
$X$ is an isotropic random vector on $\R^n$ (i.e., $\E\inr{t,X}^2 = 1$
for every $t$ in the Euclidean unit sphere),
and $Y=\inr{t_0,\cdot} + W$, for $t_0 \in \R^n$ and a symmetric random variable $W \in L_2$ that is independent of $X$ and has variance $\sigma^2$. Thus, $f^*=f_0$ and
$\F_{f^*,r} = \{\inr{t-t_0,\cdot} : \|t-t_0\|_2 \leq r\}$. If $W_1,...,W_N$ are independent copies of $W$, then by a standard symmetrization argument and since $X$ is isotropic,
\begin{align*}
& \E \sup_{u \in \F_{h,r}} \left| \frac{1}{\sqrt{N}} \sum_{i=1}^N\eps_i u(X_i) \cdot (f^*(X_i)-Y_i) \right| \leq 2\E\sup_{t \in r B_2^n} \left|\frac{1}{\sqrt{N}} \inr{\sum_{i=1}^N \eps_i W_i X_i,t} \right|
\\
= & \frac{2r}{\sqrt{N}} \E \left\|\sum_{i=1}^N \eps_i W_i X_i\right\|_2 \leq 2r \|W\|_{L_2} (\E\|X\|_2^2)^{1/2} = 2r \sigma \sqrt{n}~,
\end{align*}
where $B_2^n$ is the Euclidean unit ball in $\R^n$.
Now \eqref{eq:reg-r-n} follows from Theorem \ref{thm:main-independent}, assuming that the $L_q(\mu)$ and $L_2(\mu)$ norms are equivalent on ${\rm span}(\F)$; that is, for every $t \in \R^n$, $\|\inr{t,X}\|_{L_q} \leq L\|\inr{X,t}\|_{L_2}$. Indeed, the above shows that $\ol{r}_{\M}(c_1,f^*) \leq 2c_1^{-1}\sigma \sqrt{n/N}$. A similar argument leads to $r_E(c_1)=0$ when $N \geq cn$ for a constant $c$ that depends only on $c_1$. Also, by a volumetric estimate,
$$
\log{\cal M}(\F_{f^*,r},\eta r D) = \log {\cal M} (rB_2^n, \eta r B_2^n) \sim n \log(2/\eta)~;
$$
hence, for $N \geq cn$, $\lambda_{\Q}(c_1)=0$ and $\lambda_{\M}(c_1/\sigma,c_2) \sim \sigma \sqrt{n/N}$, implying in particular that the lower bound on the accuracy edge is $c\sigma \sqrt{n/N}$.

In other words, when $N \geq c n$, the procedure exhibits the optimal accuracy/confidence tradeoff for any $r \geq c_3\sigma \sqrt{n/N}$.
\endproof

In Section \ref{sec:examples} we present two more examples, in which we obtain the (previously unknown) optimal accuracy/confidence tradeoff. The first example studies regression in an arbitrary convex, centrally symmetric subset of $\R^n$ when the underlying measure is sub-Gaussian, but the target may be heavy-tailed. In the other example we focus on regression in $\rho B_1^n =\{t \in \R^n : \|t\|_1 \leq \rho\}$, but under significantly weaker assumptions on the underlying measure. Regression in $\rho B_1^n$ is of central importance in \emph{sparse recovery}, specifically, in the study of the \emph{basis pursuit} procedure and the \emph{LASSO}. We refer the reader to the books \cite{Kolt08,FoRa13,vdG16} for more information on sparse recovery and on these procedures.

\section{The median-of-means tournament} \label{sec:tournament}
The key to obtaining sharp estimates for both the accuracy and the confidence is identifying a procedure that is not sensitive to atypical values that may occur on a small part of the given sample. Thus, a natural starting point is the conceptually simple and attractive mean estimator, the so-called \emph{median-of-means} estimator. It was proposed, independently, by Nemirovsky and Yudin \cite{NeYu83}, Jerrum, Valiant, and Vazirani \cite{JeVaVa86},
Alon, Matias, and Szegedy \cite{AMS02}, and is defined as follows.

Let $Z_1,\ldots,Z_N$ be
independent, identically distributed real random variables with a finite second moment.
The median-of-means estimator of $\mu=\EXP Z_1$ has parameter $\delta\in [e^{1-N/2},1)$.
Setting $n=\left\lceil \ln(1/\delta)\right\rceil$, one may partition $\{1,...,N\}$ into $n$ \emph{blocks} $I_1,\ldots,I_n$, each of cardinality $|I_j|\geq \lfloor N/n\rfloor\geq 2$. Compute the sample mean in each block
$$
W_j=\frac{1}{|I_j|}\sum_{i\in I_j}Z_i
$$
and define $\wh{\mu}_N^{(\delta)}$ as the median of $W_1,\ldots,W_n$.
(If the median is not uniquely defined, here, and in the rest of the
paper, we choose the smallest one. Any other choice would work equally well.)
 It is straightforward to verify that for any $N\geq 4$,
\begin{equation}
\label{eq:mom}
\PROB\left\{|\wh{\mu}_N^{(\delta)} - \mu|>2e\sqrt{2 \var(Z)}
    \sqrt{\frac{(1+\ln(1/\delta))}{N}} \right\} \leq \delta~,
\end{equation}
where $\var(Z)$ denotes the variance of $Z$. In other words, the median-of-means estimator achieves a high (sub-Gaussian) confidence under the minimal assumption that the variance $\var(Z)$ is finite. Note that the high confidence is valid even though $Z$ can be heavy-tailed, and thus a nontrivial part of the sample $(Z_i)_{i=1}^N$ may be atypical in the sense that $Z_i$  is `far away' from $\mu$.

For properties, applications, and extensions of the
median-of-means estimator, we refer to
Bubeck, Cesa-Bianchi, and Lugosi \cite{BuCeLu13}.
Devroye, Lerasle, Lugosi, and Oliveira \cite{DeLeLuOl16}
Hsu and Sabato \cite{HsuSa13},
Lerasle and Oliveira \cite{LeOl12},
Minsker \cite{Min15},
Audibert and Catoni \cite{AuCa11}.

It is tempting to try to replace the empirical means $(1/N) \sum_{i=1}^N (f(X_i)-Y_i)^2$ by some median-of-means estimate of the risk for each $f\in \F$, and select a function in $\F$
that minimizes the estimate. However, due to the nonlinear nature of the median-of-means,
it is difficult to control the process of the estimated losses. Instead, the alternative we propose is to estimate the \emph{difference} of the risk for all pairs $f,h$ and organize a two-stage ``tournament''.

We mention here that Brownlees, Joly, and Lugosi \cite{BrJoLu15}
propose empirical minimization based on a different robust mean
estimator, a carefully designed M-estimator proposed by Catoni
\cite{Cat10}. Under general loss functions they derive analogs of Dudley's chaining bound
for the excess risk. However, the derived bounds are far from giving the optimal rate of convergence
under the squared loss.

Without loss of generality and for convenience in the notation, we use
$3N$ instead of $N$ for the sample size and assume that
$\D_{3N}=(X_i,Y_i)_{i=1}^{3N}$ is the given sample.  The sample is
split into three equal parts $(X_i,Y_i) _{i=1}^{N}$,
$(X_i,Y_i)_{i=N+1}^{2N}$ and $(X_i,Y_i)_{i=2N+1}^{3N}$.  Given the wanted degree of accuracy $r \geq r^*$, the first part of the sample--in fact, just $(X_i)_{i=1}^{N}$--is used to
estimate pairwise distances within $\F$, in a sense that will be
clarified below. The second part is used in the preliminary round of
the tournament. We show that the outcome of the preliminary round is a
set $H \subset \F$ that contains $f^*$ and possibly some other
functions whose $L_2(\mu)$ distance to $f^*$ is at most $cr$. The
final part of the tournament is a `champions league'
round. Participants in that final round are the elements in $H$ (the
`qualifiers' of the preliminary round), and the goal of that round is
to identify a function $\wh{f} \in H$ whose predictive capabilities
are almost optimal, in the sense that
$$
\E \bigl((\wh{f}(X)-Y)^2 | \D_{3N}\bigr) \leq \E(f^*(X)-Y)^2 + cr^2~,
$$
as required. This last round is only needed to guarantee the desired
excess risk. The first two rounds suffice to output a function whose
$L_2(\mu)$ distance to $f^*$ is at most $cr$: one may simply select
an arbitrary element of $H$. Also note that in the setup of Theorem
\ref{thm:main-independent}, the third round is not required as the
conditions of the theorem imply that $f^*(X)=\E(Y|X)$, and the excess
risk equals the mean squared error; thus, any $\wh{f}\in H$ has the
desired performance.

Let us now describe the three stages of the median-of-means tournament is detail.

\subsection{The `referee': the distance oracle} \label{sec-distances-oracle}
Like all good tournaments, ours too requires a `referee', whose role
is to decide whether a \emph{match} (described below) is allowed to
take place. The referee's decision is based on a \emph{distance
  oracle}--a data dependent functional that allows one to crudely
identify distances between functions in $\F$. The functional is
constructed via the median-of-means philosophy. Without loss of
generality and for ease of exposition, we may assume that $N$ is an
integer multiple of $\ell$ in the next definition.

\begin{Definition} \label{def:Med-of-means}
Let $1 \leq \ell \leq N$ and set $I^\prime_j$ to be the partition of $\{1,...,N\}$ to disjoint intervals of cardinality $\ell$. Set $k=N/\ell$ and for $v \in \R^N$ let
${\rm Med}_\ell(v)$ be the median of the means $(\ell^{-1}\sum_{i \in I^\prime_j} v_i)_{j=1}^{k}$.
\end{Definition}

Recall that one of our assumptions is an $L_q$-$L_2$ norm equivalence,
that is, that there are $q>2$ and $L \geq 1$, such that, for every $f
\in {\rm span}(\F)$, $\|f\|_{L_q} \leq L \|f\|_{L_2}$.
(In Theorem \ref{thm:main-general} we only consider $q=4$.)
 Let $\ell=\ell(q,L)$ to be specified later. For $\C_N=(X_i)_{i=1}^N$ and every $f,h \in \F$, set $v=(|f(X_i)-h(X_i)|)_{i=1}^N$ and put
$$
\Phi_{\C_N} (f,h) = {\rm Med}_\ell(v)~.
$$
The functional $\Phi$ allows one to identify distances in $\F$ in a crude (isomorphic) way, as the next theorem shows:
\begin{Proposition} \label{thm:distance-functional} 
There exist constants $\kappa,\eta,\ell, c>0$ and $0<\alpha < 1 <
\beta$, all of them depending only on $q$ and $L$ for which the
following holds. For a fixed $f^* \in \F$, let
$d^*=\max\{\lambda_{\Q}(\kappa,\eta,f^*),r_E(\kappa,f^*)\}$.
For any $r\ge d^*$,
with probability at least $1-2\exp(-cN)$, for every $f \in \F$, one has
\begin{description}
\item{$\bullet$} If $\Phi_{\C_N}(f,f^*) \geq \beta r$ then $\beta^{-1} \Phi_{\C_N}(f,f^*) \leq \|f-f^*\|_{L_2} \leq \alpha^{-1} \Phi_{\C_N}(f,f^*)$.
\item{$\bullet$} If $\Phi_{\C_N}(f,f^*) < \beta r$ then $\|f-f^*\|_{L_2} \leq (\beta/\alpha)r$.
\end{description}
\end{Proposition}

\remark
Replacing $d^*$ by the larger
$\max\{\lambda_{\Q}(\kappa,\eta),r_E(\kappa)\}$, which is independent
of $f^*$, a similar assertion to Proposition
\ref{thm:distance-functional} holds for all the pairs $f,h \in
\F$. The probability bound in that case is essentially unchanged:
$1-2\exp(-c^\prime N)$. However, Proposition
\ref{thm:distance-functional} is sufficient for our purposes.

Proposition \ref{thm:distance-functional} is an immediate modification
of Theorem 3.3 from \cite{Men16a}. For the sake of completeness we
outline the main components of its proof in the appendix.

Next we introduce the ``distance oracle'', denoted by ${\cal DO}$.
Recall the definition of $r^*$ from \eqref{eq:r*} and note that for the right choice of constants, $r^*
\geq d^*$. The distance oracle is adapted
to the wanted degree of accuracy, that is, to any fixed $r \geq 2r^*$.
\begin{Definition} \label{def:distance-oracle}
Fix $r \geq 2r^*$. Using the notation of Proposition \ref{thm:distance-functional},
if $\Phi_{\C_N}(f,h) \geq \beta r$ set ${\cal DO}(f,h)=1$, otherwise
set ${\cal DO}(f,h)=0$.
\end{Definition}
The distance oracle determines if a match between $f$ and $h$ takes
place: it does if ${\cal DO}(f,h)=1$ and it is abandoned if ${\cal
  DO}(f,h)=0$. Note that Proposition \ref{thm:distance-functional}
only shows that ${\cal DO}$ is a realistic indication of the distance
between pairs when one of the functions is the designated function
$f^*$. This serves our purposes since the designated function we are interested
in is the minimizer of the true risk in $\F$, and the success of the
procedure only requires having accurate information on matches that
involve $f^*$, even if we do not know which matches those are.

It follows from Proposition \ref{thm:distance-functional} that with
probability at least $1-2\exp(-cN)$ relative to $(X_i)_{i=1}^N$, if a
match between $f^*$ and $f$ is allowed to proceed then
$\|f-f^*\|_{L_2} \geq r$, while if it is abandoned then
$\|f-f^*\|_{L_2} \leq (\beta/\alpha)r$.

\subsection{The preliminary round} \label{sec:pre-round}

The goal of the preliminary round is to produce a subset $H \subset
\F$ that, with overwhelming probability over
the samples $(X_i,Y_i)_{i=N+1}^{2N}$, contains $f^*$ and $\|h-f^*\|_{L_2} \leq
(\beta/\alpha)r$ for any $h \in H$.

The round consists of `matches' between every pair $f,h \in \F$, and a
match can have three possible outcomes: a win by either side, or a
draw (the latter includes abandoned matches because of the ruling of
the distance oracle).

Each match is `played' using the second part of the sample,
$(X_i,Y_i)_{i=N+1}^{2N}$. The sub-sample is partitioned to $n$
 blocks $(I_j)_{j=1}^n$ of cardinality $m=N/n$ each, for a
choice of $n$ specified later. Let us note that $n$ depends on the desired degree of accuracy $r$.

\begin{description}
\item{$\bullet$} A match between $f$ and $h$ takes place if the distance oracle, using the first part of the sample $(X_i)_{i=1}^N$, declares that ${\cal DO}(f,h)=1$; otherwise, the match is abandoned and results in a draw.
\item{$\bullet$}
Each match is decided according to the $n$ blocks
  generated by the partition of $(X_i,Y_i)_{i=N+1}^{2N}$, with the
  $j$-th block played on the coordinate block $I_j$. Put
$$
B_{f,h}(j) = \frac{1}{m} \sum_{i \in I_j} \left( (f(X_i)-Y_i)^2-(h(X_i)-Y_i)^2\right), \ \ \ 1 \leq j \leq n~.
$$
The function $h$ defeats $f$ on the $j$-th block if $B_{f,h}(j) < 0$, and $f$ defeats $h$ if $B_{f,h}(j) > 0$.
\item{$\bullet$} A winner of more than $n/2$ blocks is the winner of the match. If neither function wins more than half of the blocks, the match is drawn.
\end{description}
\begin{Definition} \label{def:qualifier}
A function $f \in \F$ qualifies from the preliminary round if it has
not lost a single match; that is, it has won or drawn all its
matches. The set of ``champions'' $H$ consists of all functions qualified from the preliminary round.
\end{Definition}

The key fact regarding the outcome of the preliminary round is as follows:
\begin{Proposition} \label{thm:preliminary-round}
Under the assumptions of Theorem \ref{thm:main-general} or of Theorem \ref{thm:main-independent}, and using their notation, with probability at least
$$
1-2\exp\bigl(-c_0N \min\{1,\sigma^{-2}r^2\bigr)~,
$$
with respect to $(X_i,Y_i)_{i=1}^{2N}$, for all $h\in \F$, if ${\cal DO}(f^*,h)=1$
then $f^*$ defeats $h$. In particular, $f^* \in H$ and
for any $h \in H$, ${\cal DO}(f^*,h)=0$, and therefore $\|h-f^*\|_{L_2}
\leq (\beta/\alpha)r$.
\end{Proposition}

The proof of Proposition \ref{thm:preliminary-round} is presented in
Section \ref{sec:pre-round-proof}.
Note that Propositions \ref{thm:distance-functional} and
\ref{thm:preliminary-round} imply Theorem \ref{thm:main-independent}.
In order to prove the general result of Theorem
\ref{thm:main-general}, another round of matches is necessary to
choose a function from $H$ with small excess risk.

\subsection{Champions league} \label{sec:winners}

The goal of the second round of the tournament is to choose, among the
``champions'' selected in the preliminary round, a function with a
small excess risk. This round consists of different kind of
matches, played between functions in $H$. Since this round consists of
matches between functions in $H$, conditioned on the `good event'
from the preliminary round, every qualifier satisfies that
$\|h-f^*\|_{L_2} \leq (\beta/\alpha)r$.

The modified matches are decided using the third part of the sample
$(X_i,Y_i)_{i=2N+1}^{3N}$.
The aim is to produce a function $\wh{f} \in H$ that has a good
excess risk, namely,
$$
\E\bigl((\wh{f}(X)-Y)^2|(X_i,Y_i)_{i=2N+1}^{3N}\bigr) \leq \E(f^*(X)-Y)^2 + r_1^2~,
$$
for some $r_1$ that is bounded by a constant multiple of $r$.

Setting $\Psi_{h,f} = (h (X)-f (X))(f(X)-Y)$, the significant
observation here is that if $\E \Psi_{f^*,f}$ is not very negative,
then the prediction error associated with $f$ is small:
\begin{Lemma} \label{lemma:good-f-oracle}
For $\gamma>0$, if $f \in \F$ satisfies that $\E \Psi_{f^*,f} \geq -\gamma t^2$, then
$$
\E (f(X)-Y)^2 - \E (f^*(X)-Y)^2 \leq 2\gamma t^2~.
$$
\end{Lemma}

\proof
Observe that for every $f \in \F$
$$
(f(X)-Y)^2 - (f^*(X)-Y)^2 = (f(X)-f^*(X))^2 + 2 (f(X)-f^*(X))(f^*(X)-Y)~,
$$
and
\begin{align*}
(f(X)-f^*(X))(f^*(X)-Y) = & (f(X)-f^*(X)) \left( \left(f^*(X)-f(X)\right) + \left(f(X)-Y\right) \right)
\\
= & - (f(X)-f^*(X))^2 -\Psi_{f^*,f}~.
\end{align*}
Therefore,
$$
\E (f(X)-Y)^2 - \E(f^*(X)-Y)^2 \leq -2\E \Psi_{f^*,f} \leq 2\gamma t^2~.
$$
\endproof

The role of the ``champions league" round is to use the third
part of the sample to select $\wh{f} \in H$ for which $\E
\bigl(\Psi_{f^*,\wh{f}}|(X_i,Y_i)_{i=1}^{2N}\bigr) \geq -\gamma r_1^2$,
for a suitable constant $\gamma>0$ and $r_1$ that is proportional to
$r$.

\vskip0.4cm
The matches in the champions league consist of ``home-and-away" legs:
\begin{Definition} \label{def:champ-round} Given a sample
  $(X_i,Y_i)_{i=2N+1}^{3N}$, let $(I_j)_{j=1}^n$ be the partition of
  $\{2N+1,3N\}$ to $n$ blocks, for the same value of $n$ as in the
  preliminary round.  Let $\beta$ and $\alpha$ be as in Proposition
  \ref{thm:preliminary-round} and set $r_1=2(\beta/\alpha)r$. The
  function $f$ wins its home match against $h$ if
$$
\frac{2}{m} \sum_{i \in I_j} \Psi_{h,f}(X_i,Y_i) \geq -r_1^2/10
$$
on more than $n/2$ of the blocks $I_j$.

We select as $\wh{f}$ any ``champion" in $H$ that wins all of its home matches.
\end{Definition}

The main result regarding the champions league is as follows:
\begin{Proposition} \label{thm:winners}
Let $H \subset \F$ as above. Under the assumptions of Theorem
\ref{thm:main-general} and using its notation, with probability at least
$$
1-2\exp\bigl(-c_0N \min\{1,\sigma^{-2}r^2\}\bigr)
$$
with respect to $(X_i,Y_i)_{i=2N+1}^{3N}$, one has:
\begin{description}
\item{$\bullet$} $f^*$ wins all of its home matches, and
\item{$\bullet$} if $\E \Psi_{f^*,f} \leq -2r_1^2$, then $f$ loses its home match against $f^*$.
\end{description}
Thus, on this event, the set of possible champions is nonempty (since
it contains $f^*$), and any other champion satisfies that $\E
\Psi_{f^*,f} \geq -2r_1^2$ and therefore, by Lemma \ref{lemma:good-f-oracle},
$$
\E \bigl((\wh{f}(X)-Y)^2|(X_i,Y_i)_{i=2N+1}^{3N}\bigr) - \E(f^*(X)-Y)^2 \leq 4r_1^2~.
$$
\end{Proposition}
The proof of Proposition \ref{thm:winners} is presented in Section \ref{sec:champions-proof}.

The combination of Propositions \ref{thm:distance-functional},
\ref{thm:preliminary-round}, and \ref{thm:winners} yields the proof of
Theorem \ref{thm:main-general}.

\section{Examples} \label{sec:examples}

Before turning to the proofs of our main results, let us present some explicit examples of applications of Theorem
\ref{thm:main-general}. 

It is an unrealistic hope to obtain a simple
characterization of all involved complexity parameters in every
example.
Indeed, to get good bounds for
$\lambda_{\Q}$ and $\lambda_{\M}$, one has to obtain sharp estimates on
covering numbers of what is almost an arbitrary set and with respect
to the $L_2(\mu)$ norm for an arbitrary probability
measure $\mu$. On the other hand, $\ol{r}_{\M}$ and $r_{E}$ depend on the oscillation of
general multiplier and empirical processes, respectively. Both types of estimates
are of central importance in modern mathematics and have been the
subject of thorough research, but they are by no means completely
understood.

Having said this, there are many interesting cases in which sharp
estimates may be derived. In what follows we focus on two such
examples. The first is rather general: linear regression performed in
a convex, centrally symmetric set $T \subset \R^n$, under the
assumption that the underlying random vector $X$ is $L$-sub-Gaussian
(see the definition below). However, and unlike the results in
\cite{LeMe16}, the `noise' $Y-f^*(X)$ may be heavy-tailed.

The second example is similar to the first one and is motivated by
questions in \emph{sparse recovery}: linear regression in the set
$T=\rho B_1^n = \{t : \|t\|_1 \leq \rho\}$. The difference between
this example and the first one lies in the assumption on $X$. In the
second example $X$ is not assumed to be $L$-sub-Gaussian, but rather satisfies a much weaker moment
condition, the same condition that is needed to ensure that the \emph{basis
  pursuit} algorithm has a unique solution with the optimal number of
measurements (see \cite{LeMe16a}).

Both examples lead to explicit estimates on the accuracy and
confidence of the median-of-means tournament. The estimates are better
than the known bounds and hold with optimal confidence for accuracy
larger than $cr^*$ (though in general, the true identity of the accuracy edge is an open question). Moreover, in the second example, of $\rho B_1^n$, one may show that $r^*$ is proportional to the accuracy edge, and the optimal tradeoff holds all the way down to that value.

\subsection{Coverings and Gaussian processes}

Let $X$ be an isotropic random vector in $\R^n$, that is, for every
$t \in \R^n$, $\E\inr{X,t}^2=\|t\|_2^2$. The assumption that $X$ is isotropic
only serves clarity of the illustration. Indeed, the $L_2(\mu)$ metric endowed on $\R^n$ via the
identification of $t \in \R^n$ with the linear functional
$\inr{\cdot,t}$ is the standard Euclidean metric.  Thus, $D$ (the unit ball in $L_2(\mu)$) can be
identified with the Euclidean unit ball in $\R^n$. While it is possible to
extend the results presented below to $X$ with a general
covariance structure (in which case, $D$ is identified with an
ellipsoid in $\R^n$), the isotropic example is interesting enough to
serve as a proof of concept.

Thanks to the isotropicity assumption, if $T$ is a convex and
centrally-symmetric set
and $F=\{\inr{\cdot,t} : t \in T\}$, then for every $h \in F$,
\begin{equation*}
\lambda_{\Q}(\kappa,\eta,h) \leq  \inf\{ r: \log {\cal M}(2T \cap rB_2^n ,\eta r B_2^n) \leq \kappa^2 N\}~,
\end{equation*}
and
\begin{equation*}
\lambda_{\M}(\kappa,\eta,h) \leq  \inf\{ r: \log {\cal M}(2T \cap rB_2^n,\eta r B_2^n) \leq \kappa^2 N r^2\}~.
\end{equation*}
A standard, though sometimes suboptimal, method to estimate
covering/packing numbers relies on the theory of Gaussian processes,
specifically, on Sudakov's inequality. We formulate it
only in the case we need here. Denote by
$G=(g_i)_{i=1}^n$ a standard Gaussian vector in $\R^n$. For $T
\subset \R^n$ let
$$
\ell_*(T) = \E \sup_{t \in T} \inr{G,t}
$$
the \emph{mean-width} of $T$ with respect to the Gaussian measure.

\begin{Proposition} \label{thm:Sudakov}
(Sudakov \cite{Sud69}.)
There exists an absolute constant $c$, such that, for any $T \subset \R^n$ and every $\eps>0$,
\begin{equation} \label{eq:sudakov}
\eps \sqrt{\log{{\cal M}(T,\eps B_2^n)}} \leq c \ell_*(T)~.
\end{equation}
\end{Proposition}
Applying Proposition \ref{thm:Sudakov} it follows that
$$
\log {\cal M}(2T \cap rB_2^n ,\eta r B_2^n) \leq \left(c\frac{\ell_*(2T \cap rB_2^n)}{\eta r}\right)^2~.
$$
Hence,
\begin{equation} \label{eq:cond-gauss-lambda-Q}
\lambda_{\Q}(\kappa,\eta,h) \leq  \inf\{ r:  \ell_*(2T \cap rB_2^n) \leq (\kappa \eta/c) r \sqrt{N}\}~,
\end{equation}
and
\begin{equation} \label{eq:cond-gauss-lambda-M}
\lambda_{\M}(\kappa,\eta,h) \leq  \inf\{ r:  \ell_*(2T \cap rB_2^n) \leq (\kappa \eta/c) r^2 \sqrt{N}\}~.
\end{equation}

We emphasize again that replacing $\lambda_{\Q}$ and $\lambda_{\M}$
with these upper estimates, is, at times, suboptimal.
We refer the reader to \cite{LeMe16} for more details on this issue.

The other two parameters involved in Theorem \ref{thm:main-general},
namely, $r_E$ and $\ol{r}_{\M}$, measure the oscillation of multiplier
and empirical processes. The analysis of such processes is highly
nontrivial--even when just considering their limits as the sample
size $N$ tends to infinity, and one expects convergence to the
limiting Gaussian process (see, for example, the book \cite{Dud:book} for a
detailed exposition of such limit theorems). Because the estimates we require
are non-asymptotic, in general they are
much harder to obtain. 

The following notion makes the task of
obtaining such bounds more manageable, though still nontrivial.

\begin{Definition} \label{def:subgaussian-vector}
A random vector $X$ in $\R^n$ is $L$-sub-Gaussian if for every $t \in \R^n$ and any $p \geq 2$,
$$
\|\inr{X,t}\|_{L_p} \leq L \sqrt{p} \|\inr{X,t}\|_{L_2}~.
$$
\end{Definition}
Note that if $X$ is $L$-sub-Gaussian and isotropic, then for every $t
\in \R^n$, $\|\inr{X,t}\|_{L_p} \leq L\sqrt{p} \|t\|_2$, because $\|\inr{X,t}\|_{L_2}=\|t\|_2$.

The simplest examples of isotropic, $L$-sub-Gaussian random vectors
are vectors with independent, mean-zero, variance $1$ components that
are sub-Gaussian.  For instance, the standard Gaussian vector
$(g_i)_{i=1}^n$, and $(\eps_i)_{i=1}^n$, whose components are
independent, symmetric random signs are both $L$-sub-Gaussian for a
constant that is independent of the dimension $n$. Another family
of examples consists of the random vectors whose
density is uniform on sets of the form $\{t : \|t\|_p \leq cn^{1/p}\}$
for some $p \geq 2$, normalized to have volume $1$. Again, $L$ is an absolute
constant, independent of $n$ and $p$ (see, e.g., \cite{BaKo03,BGMN05}).

The reason for considering a sub-Gaussian random vector $X$ is that,
by Talagrand's theory of \emph{generic chaining}
(see the book \cite{Tal14} for an extensive exposition on the subject), the oscillations in question may be controlled using the oscillation of the corresponding Gaussian process. For example, the
next result describes how the expected supremum of a multiplier
process is upper bounded in terms of the Gaussian mean-width
$\ell_*$. Although it is formulated in $\R^n$, it holds in a
far more general context (see \cite{Men16b}).

\begin{Proposition} \label{thm:osc} 
Let $q >2$ and let $X$ be an isotropic, $L$-sub-Gaussian random vector in
  $\R^n$. There exists a constant $c=c(q)$ such that the following
  holds. Let $\xi \in L_q$ be a random variable (not necessarily
  independent of $X$) and let $(X_i,\xi_i)_{i=1}^N$ be independent
  copies of $(X,\xi)$. Then, for any $T  \subset \R^n$,
\begin{equation} \label{eq:multi}
\E \sup_{t \in T}  \left|\frac{1}{\sqrt{N}} \sum_{i=1}^N \eps_i \xi_i \inr{X_i,t} \right| \leq c L \|\xi\|_{L_q} \ell_*(T)~,
\end{equation}
where $(\eps_i)_{i=1}^N$ are independent, symmetric signs that are
independent of $(X_i,\xi_i)_{i=1}^N$.
\end{Proposition}

Note that if $\xi$ is heavy-tailed there is no hope of obtaining a
high-probability version of Proposition \ref{thm:osc}.
 In fact, if all one knows is that $\xi$ belongs to $L_q$, one cannot hope that
$$
\sup_{t \in T}  \left|\frac{1}{\sqrt{N}} \sum_{i=1}^N \xi_i \inr{X_i,t} \right| \leq c L \|\xi\|_{L_q} \ell_*(T)
$$
with a probability higher than $1-c_1/N^{(q/2)-1}$.

This exhibits, yet again, the weakness of Theorem
\ref{thm:best-known}, which is based on the ``in-probability"
parameter $r_{\M}$. Setting $f^*(X)=\inr{X,t^*}$ and if $\xi=f^*(X)-Y$
is heavy-tailed, the oscillation is simply too big on an event with
the required confidence \eqref{eq:optimal-confidence}, but the mean
oscillation is well behaved.

Applying Proposition \ref{thm:osc}, it is evident that
\begin{equation} \label{eq:r-E-subgaussian}
r_{E}(\kappa) \leq \inf\left\{ r: \ell_*(2T \cap rB_2^n) \leq (\kappa/cL) \sqrt{N} r\right\}
\end{equation}
and
\begin{equation}
\ol{r}_{\M}(\kappa,t) = \inf\left\{ r: \|Y-\inr{X,t}\|_{L_q} \ell_*(2T \cap rB_2^n) \leq (\kappa/cL) \sqrt{N} r^2\right\}~.
\end{equation}
Thus, assuming that $\|Y-\inr{X,t}\|_{L_q} \leq L \|Y-\inr{X,t}\|_{L_2}$, as we do in Assumption \ref{ass:moment}, it follows that
$$
{\wt r}_{\M}(\kappa,\sigma) \leq  \inf\left\{ r: \ell_*(2T \cap rB_2^n) \leq (\kappa/cL^2 \sigma) \sqrt{N} r\right\}~.
$$
\begin{Definition}
For constants $c_1$, $c_2$ and $\sigma$, let
$$
s_{\M}(c_1,\sigma)=\inf\left\{ r: \ell_*(2T \cap rB_2^n) \leq (c_1/\sigma) \sqrt{N} r\right\}~,
$$
and
$$
s_{\Q}(c_2)=\inf\left\{ r: \ell_*(2T \cap rB_2^n) \leq c_2 \sqrt{N} r\right\}~.
$$
\end{Definition}
Hence, for every $t^* \in T$,
\[
\max\{\lambda_{\Q}(c_1,c_2,t^*),\lambda_{\M}(c_1/\sigma,c_2,t^*),r_E(c_1,t^*),\ol{r}_{\M}(c_1,t^*) \} \leq \max\{s_{\Q}(c_3,\sigma),s_{\M}(c_4/\sigma)\}~,
\]
for constants $c_1,c_2,c_3,c_4$ that depend only on $q$ and $L$.

We obtain the following consequence of Theorem
\ref{thm:main-general}:

\begin{Theorem} \label{thm:subgaussian-X-heavy-noise}
Let $X$ be an isotropic, $L$-sub-Gaussian random vector in $\R^n$, and
let $T \subset \R^n$ be a convex, centrally-symmetric set. Let $q>2$
and assume that $\|Y-\inr{X,t}\|_{L_q} \leq L\|Y-\inr{X,t}\|_{L_2}$
for every $t \in T$. Then, with probability at least $1-2\exp(-c_1 N
\min\{1,\sigma^{-2}s_{M}^2(c_2,\sigma)\})$, the median-of-means
tournament produces $\wh{t}$ such that
\[
\|\wh{t}-t^*\|_{L_2} \leq c_3 s^* \ \ {\rm  and \ \ }
\E \bigl((\inr{X,\wh{t}}-Y)^2 |\D_N \bigr) \leq \E (\inr{X,t^*}-Y)^2 + (c_4 s^*)^2~,
\]
where
\[
s^*=\max\{s_{\M}(c_2,\sigma),s_{\Q}(c_5)\}~,
\]
and the constants $c_1,...,c_5$ depend only on $L$ and $q$.

Moreover, for any $s \geq s^*$,
\[
\|\wh{t}-t^*\|_{L_2} \leq c_3 s \ \ {\rm  and \ \ }
\E \bigl((\inr{X,\wh{t}}-Y)^2 |\D_N \bigr) \leq \E (\inr{X,t^*}-Y)^2 + (c_4 s)^2~,
\]
with probability at least $1-2\exp(-c_1 N \min\{1,\sigma^{-2}s^2\})$, exhibiting the optimal accuracy-confidence tradeoff.
\end{Theorem}
Theorem \ref{thm:subgaussian-X-heavy-noise} improves Theorem A from
\cite{LeMe16}, where it was shown that ERM produces $\wh{t}$ with the same
accuracy and confidence as in the first part of Theorem
\ref{thm:subgaussian-X-heavy-noise}, but only when
$\|Y-\inr{X,t}\|_{L_q} \leq L\sqrt{q}\|Y-\inr{X,t}\|_{L_2}$ \emph{for
  every} $q>2$ and every $t \in T$. In other words, Theorem A from
\cite{LeMe16} is based on the assumption that each $Y-\inr{X,t}$ is an
$L$-sub-Gaussian random variable, and holds only for the accuracy level $s^*$. In contrast, Theorem
\ref{thm:subgaussian-X-heavy-noise} shows that the median-of-means
tournament
performs in an optimal way in heavy-tailed situations that are totally out of reach for ERM and for the entire range $s \geq cs^*$.

Observe that the only range of accuracies in which Theorem \ref{thm:subgaussian-X-heavy-noise} is (perhaps) suboptimal, is when
\begin{equation} \label{eq:subopt-range}
\lambda^*=\max\{\lambda_{\Q}(\kappa_1,\eta_1),\lambda_{\M}(\kappa_2/\sigma,\eta_2)\} \leq s \leq s^*
\end{equation}
for well chosen values of $\kappa_i,\eta_i$, $i=1,2$; that is, for values that are larger than the known lower estimate on the accuracy edge for such problems. As noted in \cite{LeMe16}, there are many examples in which
$\lambda^*$ and $s^*$ are equivalent (roughly speaking, this happens when
Sudakov's inequality is sharp). In such cases the median-of-means tournament is optimal in the entire range of accessible accuracies.

One important class of sets in which this equivalence is true is $\rho
B_1^n = \{t : \|t\|_1 \leq \rho\}$ (see \cite{LeMe16} for the proof).
In light of Theorem \ref{thm:subgaussian-X-heavy-noise}, the
median-of-means tournament performs in an optimal way in $\rho
B_1^n$. Moreover, it turns out that one may relax the sub-Gaussian
assumption on $X$ and still obtain the optimal behaviour $\rho
B_1^n$, as we show next.

\subsection{$\rho B_1^n$ -- Sparse recovery sets}
It is well understood that classical sparse recovery procedures, such as \emph{basis pursuit} or \emph{LASSO} relay heavily on the geometry of $B_1^n$. Indeed, LASSO selects $\wh{t}$, the minimizer in $\R^n$ of the functional
$$
t \to \frac{1}{N} \sum_{i=1}^N \left(Y_i-\inr{X_i,t}\right)^2 + \lambda \|t\|_1~.
$$
Being a penalized version of ERM, the
analysis of LASSO is equivalent to the study of ERM in the sets $\rho B_1^n$ for an arbitrary choice of $\rho$.

While we defer the question of a ``LASSO-tournament" procedure to
future work, it is clear that the first step in that direction is to
explore the median-of-means tournament in $\rho B_1^n$. Instead of the
sub-Gaussian assumption used in Theorem
\ref{thm:subgaussian-X-heavy-noise}, the assumption we use follows the
path of \cite{LeMe16a}:
\begin{Assumption} \label{ass:CS-moments}
Let $X$ be an isotropic random vector and $\kappa \geq 1$. Assume that
for every $t \in \R^n$ and any $2 \leq p \leq \kappa \log n$,
$\|\inr{X,t}\|_{L_p} \leq L \sqrt{p} \|\inr{X,t}\|_{L_2}$.
\end{Assumption}
Note that the coordinates of $X$ need not be independent and $X$ may
be far from being an $L$-sub-Gaussian random vector. Indeed, it is
required that linear functionals $\inr{\cdot,t}$ satisfy a
sub-Gaussian moment growth only up to the logarithm of the dimension,
and it is possible that some do not have any higher moment beyond $p=\kappa
\log n$.

It turns out (see \cite{LeMe16a}) that if $X$ satisfies Assumption
\ref{ass:CS-moments}, then
$N^{-1/2} \sum_{i=1}^N \inr{X_i,\cdot}e_i$,
the random matrix whose rows are independent copies of $X$, exhibits
the best possible sparse recovery features. For example, one requires
$N \sim s\log(en/s)$ random measurements $\inr{X_i,t}$ to recover any
$s$-sparse vector $t$ using basis pursuit, and a similar type of
estimate holds in the ``noisy" setup for the LASSO. Moreover, one
cannot
relax the moment condition in Assumption \ref{ass:CS-moments} and
still get the same recovery properties.

Of course, since both basis pursuit and LASSO are variations of ERM,
they suffer from the same weaknesses as ERM. As such, when the given
measurements are $(\inr{X_i,t})_{i=1}^N$ for a heavy-tailed $X$, the
confidence with which the recovery properties hold is suboptimal, and
very different from the confidence one has when $X$ is the standard
Gaussian vector in $\R^n$.

We show that as far as regression in $\rho B_1^n$ goes, the
median-of-means tournament
yields the optimal, ``Gaussian" behaviour even when $X$ only
satisfies Assumption \ref{ass:CS-moments} and $Y-\inr{t^*,X}$ is
heavy-tailed. To this end, we need sharp bounds on the
parameters that are used to define $r^*$ in Theorem
\ref{thm:main-general}.

As was noted earlier, $\lambda_{\Q}$ and $\lambda_{\M}$ depend only on
the covariance structure endowed on $\R^n$ by $L_2(\mu)$, and since
$X$ is isotropic, the $L_2(\mu)$ metric corresponds to the standard
Euclidean norm. Therefore, the difficulty lies in bounding $r_E$ and
$\ol{r}_{\M}$, and specifically, in extending Proposition
\ref{thm:osc} beyond the $L$-sub-Gaussian case. While Theorem
\ref{thm:osc} is general and holds for any subset of $\R^n$, here it
is needed for very specific sets, namely $T=\rho B_1^n \cap r B_2^n$. Such
indexing sets fall within the scope of Theorem 1.6 from \cite{Men16c}.

\begin{Proposition} \label{thm:osc-weak-moment}
Let $q >2$ and let $X$ satisfy Assumption
\ref{ass:CS-moments} for $\kappa=c_1(q)$.
Let $\xi \in L_q$ be a random variable (not necessarily
  independent of $X$) and let $(X_i,\xi_i)_{i=1}^N$ be independent
  copies of $(X,\xi)$. Then
\begin{equation} \label{eq:multi1}
\E \sup_{t \in \rho B_1^n \cap r B_2^n}  \left|\frac{1}{\sqrt{N}}\sum_{i=1}^N \eps_i \xi_i \inr{X_i,t} \right| \leq c_2(q) L \|\xi\|_{L_q} \ell_*(\rho B_1^n \cap r B_2^n)~.
\end{equation}
\end{Proposition}
Thus, Assumption \ref{ass:CS-moments} suffices to ensure that $r_E$ and ${\wt r}_{\M}(\kappa,\sigma)$ may be controlled as if $X$ were $L$-sub-Gaussian. All that remains is to estimate
$$
\log {\cal M}(\rho B_1^n \cap rB_2^n, \eta r B_2^n) \ \ {\rm and} \ \ \ell_*(\rho B_1^n \cap r B_2^n)
$$
which are well-understood quantities.

To put them in a more familiar form, set $\sqrt{s}=\rho/r$ and observe that
$$
\ell_*(\rho B_1^n \cap r B_2^n)=r \ell_*((\rho/r)B_1^n \cap B_2) = r \ell_*(\sqrt{s}B_1^n \cap B_2^n)
$$
and
$$
{\cal M}(\rho B_1^n \cap rB_2^n, \eta r) = {\cal M}((\rho/r)B_1^n \cap B_2^n, \eta B_2^n) = {\cal M}(\sqrt{s}B_1^n \cap B_2^n, \eta B_2^n)~.
$$
Recall that if $1 \leq s \leq n$ and $V_s$ is the set of $s$-sparse
vectors in the Euclidean unit sphere (i.e., those with at most $s$
nonzero components), then ${\rm conv}(V_s) \subset
\sqrt{s}B_1^n \cap B_2^n \subset C \cdot {\rm conv}(V_s)$ for a
suitable absolute constant $C$. Using what are by now standard
estimates (see, e.g., \cite{LeMe16}),
$$
\ell_*(\sqrt{s}B_1^n \cap B_2^n) \sim \sqrt{s\log(en/s)}, \ \ {\rm and} \ \ \log{\cal M}(\sqrt{s}B_1^n \cap B_2^n, \eta B_2^n) \sim s\log(en/\eta s)~.
$$
The estimates are simpler outside the range $1 \leq s \leq n$: when
$\rho/r \leq 1$ then $\rho B_1^n \cap r B_2^n = \rho B_1^n$ and when
$\rho/r \geq \sqrt{n}$ then $\rho B_1^n \cap r B_2^n = r B_2^n$.
Again, the required estimates on ${\cal M}$ and $\ell_*$ are standard
and may be found, for example, in \cite{LeMe16}.

Using these observations, and with the same (tedious) computation as
in \cite{LeMe16}, one obtains the following:
let $c_1$ and $c_2$ be well-chosen absolute constants and set
\begin{equation*}
v_{\M}^2 =
\begin{cases}
\frac{\rho \sigma}{\sqrt{N}} \sqrt{\log\left(\frac{2c_1n\sigma}{\sqrt{N}\rho}\right)} & \mbox{if } \ \  N \leq c_1n^2 \sigma^2/\rho^2
\\
\\
\frac{\sigma^2 n}{N} & \mbox{if} \ \ N >c_1n^2 \sigma^2/\rho^2~,
\end{cases}
\end{equation*}
and
\begin{equation*}
v_{\Q}^2=
\begin{cases}
\frac{\rho^2}{N} \log\left(\frac{2c_2n}{N}\right) &  \mbox{if} \ \ N \leq c_2 n~,
\\
\\
0 & \mbox{if} \ \ N > c_2 n~.
\end{cases}
\end{equation*}
Then
\begin{description}
\item{$\bullet$} $\lambda^*$, the lower estimate on the accuracy edge, satisfies $\lambda^* \geq c_3\max\{v_{\Q},v_{\M}\}$; thus, there is no learning procedure in $\rho B_1^n$ that can perform with a better accuracy than $c_3\max\{v_{\Q},v_{\M}\}$ with a higher confidence than $3/4$;
\item{$\bullet$} For any $v \geq c_4\max\{v_{\Q},v_{\M}\}$, the
  median-of-means tournament achieves the accuracy $v$ with the optimal confidence $1-2\exp(-c_5N \min\{1,\sigma^{-2}v^2\})$,
    thus exhibiting the optimal accuracy/confidence tradeoff up to a level that is proportional to $\lambda^*$.
\end{description}
Formally:
\begin{Corollary} \label{cor:B-1-n-weak moments}
Let $q>2$ and assume that $X$ satisfies Assumption \ref{ass:CS-moments} with a constant $\kappa=c(q)$. Assume further that for every $t \in \rho B_1^n$, $\|Y-\inr{X,t}\|_{L_q} \leq L \|Y-\inr{X,t}\|_{L_2}$.
Then for every $v \geq c_4\max\{v_{\Q},v_{\M}\}$, with probability at
least $1-2\exp\left(-c_5N \min\{1,\sigma^{-2}v^2\}\right)$, the
median-of-means tournament produces $\wh{t} \in \rho B_1^n$ that satisfies
$$
\|\wh{t}-t^*\|_{\ell_2^n} \leq c_6v \ \ {\rm and} \ \ \E \bigl((\inr{X,\wh{t}}-Y)^2 |\D_N \bigr) \leq \E (\inr{X,t^*}-Y)^2 + (c_6v)^2~,
$$
for constants $c_4,c_5,c_6$ that depend only on $q$ and $L$.
\end{Corollary}
The advantage of the median-of-means tournament over ERM is clear: it
performs in $\rho B_1^n$ with the optimal accuracy and confidence,
starting a constant factor away from the level of accuracy that can be
attained only with constant confidence, and it does so under a
heavy-tailed assumption both on $X$ and on $Y$. In contrast, ERM (which was the ``record holder" prior to this work) achieves the optimal performance only in a purely sub-Gaussian setup and does so only for one level of accuracy, of the order of $\max\{v_{\Q},v_{\M}\}$.

\section{Proofs}
Let us begin by considering the structure of the various indexing sets involved in the proofs, paying particular attention to the way their structure affects the regularity of the parameters we defined earlier.

For any class $\F$ and every $h \in \F$, $\F - h \subset \F - \F$, and in particular,
$$
\F_{r,h} = {\rm star}(\F-h,0) \cap rD \subset {\rm star}(\F-\F,0) \cap rD.
$$
Thus, for all the complexity parameters defined above, one may avoid
the need to take the supremum over all possible choices of centres by
considering a slightly larger indexing set, namely, ${\rm  star}(\F-\F,0) \cap rD$.
Moreover, if $\F$ is convex, then $\F-\F$
is both convex and centrally symmetric, and if $\F$ happens to be convex and centrally symmetric then $\F-\F=2\F$ and ${\rm star}(\F-\F,0) \cap rD = 2\F \cap rD$.

The fact that ${\rm star}(\F-h,0)$ is star-shaped around $0$ leads to important regularity properties of the parameters we use. Recall that $S$ is the unit sphere in $L_2(\mu)$ and observe that if $V$ is star-shaped around $0$ and $v \in V \cap rS$, then for every $r^\prime \leq r$, $V$ contains a `scaled-down' version of $v$, of norm $r^\prime$: for $\alpha=r^\prime/r <1$, $\alpha v \in V$. Hence, if $\phi(r) = \E \sup_{v \in V \cap rD} \bigl|\sum_{i=1}^N \eps_iv(X_i) \bigr|$ and $r^\prime \leq r$, then $\phi(r^\prime) \geq (r^\prime/r) \phi(r)$.

This argument shows that when $r>r_E(\kappa,h)$,
$$
\E \sup_{v \in \F_{h,r}} \left|\frac{1}{\sqrt{N}}\sum_{i=1}^N \eps_iv(X_i) \right| \leq \kappa \sqrt{N} r~,
$$
and when $r < r_E(\kappa,h)$, the reverse inequality holds. In a similar fashion, if $r \geq {\ol r}_{\M}(\kappa,h)$, then
$$
\E \sup_{v \in \F_{h,r}} \left|\frac{1}{\sqrt{N}} \sum_{i=1}^N \eps_iv(X_i)(h(X_i)-Y_i) \right| \leq \kappa \sqrt{N} r^2~,
$$
whereas if $r < {\ol r}_{\M}(\kappa,h)$, the reverse inequality holds.

Similar arguments apply for $\lambda_{\Q}$ and $\lambda_{\M}$: if $V$
is star-shaped around $0$ and if $\{v_1,...,v_M\}$ is an $\eta
r$-separated subset of $V \cap r D$ and $0<\alpha <1$,
then $\{\alpha v_1,...,\alpha v_M\}$ is an $\eta (\alpha r)$-separated subset of $V \cap (\alpha r)D$. Therefore,
\begin{equation} \label{eq:monotone-covering}
\log {\cal M}(V \cap \alpha rD, \eta \alpha r D) \geq \log {\cal M}(V \cap rD, \eta rD)~,
\end{equation}
implying that the function $r \to \log {\cal M}(V \cap rD, \eta rD)$ is monotone decreasing. Moreover,
if $\log {\cal M}(V \cap r^\prime D, \eta r^\prime D) \leq \kappa^2
N$, or if
$\log {\cal M}(V \cap r^\prime D, \eta r^\prime D) \leq \kappa^2
N(r^\prime)^2$,
the same is true for every $r \geq r^\prime$. Therefore, if $r > \lambda_{\Q}(\kappa,\eta,h)$ then
$$
\log {\cal M}(\F_{h,r}, \eta r D) \leq \kappa^2 N~,
$$
and if $r<\lambda_{\Q}(\kappa,\eta,h)$, the reverse inequality holds.
In a similar fashion, if $r > \lambda_{\M}(\kappa,\eta,h)$ then
$$
\log {\cal M}(\F_{h,r}, \eta r D) \leq \kappa^2 N r^2,
$$
while if $r<\lambda_{\Q}(\kappa,\eta,h)$, the reverse inequality holds.

These observations allow one to choose a level of accuracy such as $r^*$ by
``intersecting" multiple conditions like the ones appearing in the
definitions of $\lambda_{\Q},\lambda_{\M},r_E$ and $\ol{r}_{\M}$,
simply because each one of the required inequalities holds for any
level larger than $\lambda_{\Q},\lambda_{\M},r_E$ and $\ol{r}_{\M}$,
respectively.  An additional feature is that one may ``combine"
conditions by decreasing the constants involved in the
definitions. For example, $\max\{r_{E}(\kappa_1,h),r_{E}(\kappa_2,h)\} =
r_{E}(\min\{\kappa_1,\kappa_2\},h)$ and similar observations hold for
$\lambda_{\Q},\lambda_{\M}$ and $\ol{r}_{\M}$.

For the sake of transparency, we do not specify all the constants
involved in the definitions of the fixed points right from the
start. Instead, we collect conditions on these constants and use the
fact that one may ``combine" and ``intersect" them. It turns out that
all the constants depend on only two parameters: $q>2$, for which the
$L_q$ norm is equivalent to the $L_2$ norm, and the
constant $L$. In what follows, we denote by $c(q,L)$ a constant that
depends only on $q$ and $L$, for the values of $q$ and $L$ in Theorem
\ref{thm:main-general} (where only $q=4$ is considered) or in Theorem
\ref{thm:main-independent}, when $q$ can be arbitrarily close to $2$--at the price of worse constants, of course.

With this in mind, let $\kappa_1,\kappa_2, \kappa_3$ and $\eta$ be constants that will be specified later and that depend only on $q$ and $L$. Fix $f^* \in \F$ and consider $r^*$ for which, for every $r > r^*$,
\begin{equation} \label{eq:r*-emp}
\E \sup_{f \in \F_{f^*,r}} \left| \frac{1}{\sqrt{N}} \sum_{i=1}^N \eps_i (f-f^*)(X_i) \right| \leq \kappa_1 \sqrt{N}r~,
\end{equation}

\begin{equation} \label{eq:r*-multi}
\E \sup_{f \in \F_{f^*,r}} \left| \frac{1}{\sqrt{N}} \sum_{i=1}^N \eps_i (f^*(X_i)-Y_i)(f-f^*)(X_i) \right| \leq \kappa_2 \sqrt{N}r^2
\end{equation}
and
\begin{equation} \label{eq:r*-entropy}
\log{\cal M}(\F_{f^*,r},\eta rD) \leq \kappa_3^2 N \min\{1,r^2/\sigma^2\}~.
\end{equation}
The value of $r^*$ from \eqref{eq:r*} will be given by selecting right
constants $\kappa_1,\kappa_2, \kappa_3$ and $\eta$. We consider an accuracy 
$r \geq 2r^*$, and for the rest of this section we fix its value.

Let us introduce the main parameter $n$ used in the tournament,
that is, the
number of blocks into which the second part of the sample $(X_i,Y_i)_{i=N+1}^{2N}$, and the third part of the sample $(X_i,Y_i)_{i=2N+1}^{3N}$ are partitioned.
To this end, let $0<\tau<1/4$, set $0<\theta \leq \tau$, to be specified later, and put $\lambda_{\M}=\lambda_{\M}(\kappa_3/\sigma,\eta,f^*)$. Set
\begin{equation} \label{eq:number-of-blocks}
n=\theta N \min\left\{1,\left(\frac{r}{\sigma}\right)^2\right\}~.
\end{equation}

Hence, $n$ depends on the wanted accuracy $r$. Also, $n \leq \tau N $ and without loss of generality we may assume
that both $n$ and $m=N/n$ are integers.  Also, note that $m \geq 1/\tau$.

\subsection{The preliminary round--proof} \label{sec:pre-round-proof}

In this section we prove Proposition \ref{thm:preliminary-round}.

Recall that by Proposition \ref{thm:distance-functional} and the
resulting condition on $r$, one has that, with probability at least
$1-2\exp(-c_0N)$ with respect to $(X_i)_{i=1}^N$, if a match involving
$f^*$ is allowed to take place (i.e., if ${\cal DO}(f^*,f)=1$), then
$\|f^*-f\|_{L_2} \geq r$; and if the match is abandoned (i.e.,
${\cal DO}(f^*,h)=0$), then $\|f^*-f\|_{L_2} \leq
(\beta/\alpha)r$. The constants $\alpha<1<\beta$ and $c_0$ depend
only on $L$ and $q$.

Therefore, to establish Proposition \ref{thm:preliminary-round} it is
enough to show that with probability $1-2\exp(-cn)$,
if $\|f-f^*\|_{L_2} \geq r$, then $f^*$ wins its match against
$f$.
In other words, that $B_{f,f^*}(j)>0$ for the majority of the blocks in any such match.

Let us begin by exploring the situation in a match between $f^*$ and $f$, knowing that $\|f^*-f\|_{L_2} \geq r \geq 2r^*$, as above.

Clearly, for every $f,h \in \F$,
\begin{equation*}
(f(X)-Y)^2-(h(X)-Y)^2 = (f(X)-h(X))^2+2(f(X)-h(X))\cdot (h(X)-Y)~.
\end{equation*}
Set
$$
{\Q}_{f,h} = \frac{1}{m}\sum_{i=1} (f(X_i)-h(X_i))^2 \ \ {\rm and} \ \ {\M}_{f,h}=\frac{2}{m}\sum_{i=1}^m (f(X_i)-h(X_i))\cdot (h(X_i)-Y_i)~.
$$
We introduce the notation $\xi=f^*(X)-Y$ and $\xi_i=f^*(X_i)-Y_i$.
Partitioning the sample $(X_i,Y_i)_{i=N+1}^{2N}$ into the $n$ blocks $(I_j)_{j=1}^n$, each one of cardinality $m$, one has that for $1 \leq j \leq n$, ${\M}_{h,f^*}(j)=\frac{2}{m}\sum_{i \in I_j}^m \xi_i(f (X_i)-f^*(X_i))$ and
$$
B_{f,f^*}(j) = \frac{1}{m} \sum_{i \in I_j} (f (X_i)-f^*(X_i))^2 + \frac{2}{m} \sum_{i \in I_j} \xi_i (f (X_i)-f^*(X_i))={\Q}_{f,f^*}(j)+{\M}_{f,f^*}(j)~.
$$

It follows that $f^*$ defeats $f$ in the $j$-th block if $\Q_{f,f^*}(j)+\M_{f,f^*}(j) > 0$, and, because $\tau < 1/4$, $f^*$ wins the match if $\Q_{f,f^*}(j) \geq C\|f-f^*\|^2_{L_2}$ on more than $(1-\tau)n$ of the blocks, while $\M_{f,f^*}(j) \leq -C\|f-f^*\|_{L_2}^2$ on at most
$\tau n$ of the blocks, .

This is summarized in the following lemma, established in the next two
sections, and for the right choice of the constants $\kappa_1,\kappa_2,\kappa_3$,
$\eta,\tau$, and $\theta$ that depend only on $L$ and $q$.

\begin{Lemma} \label{lemma:components-of-main}
There exists an absolute constant $c$ and a constant $C_1=C_1(L,q)$
for which the following holds. With probability at least $1-2\exp(-c
\tau^2 n)$, for every $f \in \F$ with $\|f-f^*\|_{L_2} \geq r$,
$$
\left|\left\{j : {\Q}_{f,f^*}(j) \geq C_1 \|f-f^*\|_{L_2}^2 \right\} \right| \geq \left(1-\tau\right)n
$$
and
$$
\left|\left\{j : {\M}_{f,f^*}(j) \leq -\frac{3C_1}{4}\|f-f^*\|_{L_2}^2\right\} \right| \leq \tau n~.
$$
\end{Lemma}

The analysis takes place in the set
$$
\F_1={\rm star}(F,f^*) \cap (f^* + rS)~,
$$
consisting of all the functions in the star-shaped hull of $\F$ and
$f^*$ whose distance to $f^*$ is precisely $r$. Once the
estimates in Lemma \ref{lemma:components-of-main} are verified for
functions in $\F_1$, extending them to the set $\{f \in {\rm star}(\F,f^*) : \|f-f^*\|_{L_2} \geq r\}$, is straightforward, by
invoking homogeneity properties of ${\Q}_{f,f^*}$ and ${\M}_{f,f^*}$
in $f-f^*$, and because ${\rm star}(\F,f^*)$ is star-shaped around
$f^*$.

\subsubsection*{The quadratic component} \label{sec:quad}
Here we establish the first part of Lemma \ref{lemma:components-of-main},
that states that, on an event of  high probability, whenever $\|f-f^*\|_{L_2}$ is
sufficiently large, a significant majority of the
$({\Q}_{f,f^*}(j))_{j=1}^n$ are at least a large fixed proportion of
$\|f-f^*\|^2_{L_2}$. The size of the fixed proportion depends only on the
small-ball property satisfied by the class $\F$, which, in turn,
follows from the $L_q$-$L_2$ norm equivalence in ${\rm span}(F)$.
Indeed, if $\|f-f^*\|_{L_q} \leq L\|f-f^*\|_{L_2}$,
then by the Paley-Zygmund inequality (see, e.g., \cite{deGi99}),
\begin{equation} \label{eq:small-ball-in-proof}
\PROB(|f-f^*|(X) \geq \kappa_0\|f-f^*\|_{L_2}) \geq \rho_0
\end{equation}
for constants $\kappa_0$ and $\rho_0$ that depend only on $L$ and $q$.

Recall that $(I_j)_{j=1}^n$ are the $n$ blocks, each one of
cardinality $m=N/n$, that $\F_1=\{f \in {\rm star}(\F,f^*) :
\|f-f^*\|_{L_2}=r\}$, and that $m \geq 1/\tau$. For $t>0$ and a function $u$, set
\begin{equation*}
R_j(u,t)= \left| \left\{ i \in I_j : |u(X_i)| \geq t \right\} \right|
= \sum_{i \in I_j} \IND_{\{ |u(X_i)| \geq t \}}~.
\end{equation*}
Thanks to \eqref{eq:small-ball-in-proof}, it is evident that if $f \in \F_1$ and $1 \leq j \leq n$ then, with probability at least $1-2\exp(-c_0\rho_0m)$,
$$
R_j(f-f^*,\kappa_0r) \geq \frac{m \rho_0}{2}~.
$$
Recalling that $m \geq 1/\tau$, if $\tau \leq c_1(\rho_0)$ then
$$
1-2\exp(-c_0\rho_0m) \geq 1-2\exp(-c_0\rho_0/\tau) \geq 1-\tau/12~,
$$
and $(\IND_{\{R_j(f-f^*,\kappa_0 r) \geq m\rho_0/2\}})_{j=1}^n$ are
i.i.d.\ Bernoulli random variables with mean at least $1-\tau/12$.
By standard concentration properties of Binomial distributions, with probability at least $1-2\exp(-c_2\tau^2 n)$,
\begin{equation} \label{eq:quad-good-event-1}
\frac{1}{n}\sum_{j=1}^n \IND_{\{R_j(f-f^*,\kappa_0 r) \geq m\rho_0/2\}} \geq 1-\tau/6~.
\end{equation}
In other words, at least $(1-\tau/6)n$ of the blocks $I_j$
have at least $m\rho_0/2$ coordinates of size at least $\kappa_0 r$.
Moreover, in light of the high probability estimate with which
\eqref{eq:quad-good-event-1} holds,
by the union bound, the same property is satisfied uniformly by any
subset of $\F_1$ of cardinality at most $\exp(c_2\tau^2 n/2)$.
The set we consider is a maximal $\eta r$-separated subset of $\F_1$
with respect to the $L_2(\mu)$ norm
and for a well-chosen $\eta$ specified later.

Let $H_1$ be such a maximal separated set. Since
$$
\F_1 \subset f^*+({\rm star}(\F-f^*,0) \cap rD)=f^*+\F_{f^*,r}~,
$$
it follows from the translation invariance of packing numbers that
$$
|H_1| \leq {\cal M}(\F_{f^*,r},\eta r)~.
$$
Thus, to obtain a uniform control over points in $H_1$, it suffices to show that
$$
\log {\cal M}(F_{f^*,r},\eta r) \leq (c_2/2) \tau^2 n~.
$$
Recalling \eqref{eq:r*-entropy} and the choice of $n$, it suffices to verify that
$$
\kappa_3^2 N \min\{1,\sigma^{-2}r^2\} = \kappa_3^2n/\theta \leq  (c_2/2) \tau^2 n~,
$$
which is the case as long as
\[
%begin{equation} \label{eq:kappa3-cond1}
\kappa_3^2 \leq (c_2/2) \theta \tau^2~.
\]
Observe that since $H_1$ is maximal, it is also an $\eta r$-net in
$\F_1$, that is,
every $f \in \F_1$ has some $\pi f \in H_1$ for which $\|f-\pi
f\|_{L_2} \leq \eta r$.
Consider the following event:
\begin{description}
\item{$(A)$} \eqref{eq:quad-good-event-1} holds for every $f \in H_1$, and
\item{$(B)$} $\sup_{f \in \F_1} (1/n)\sum_{i =1}^n \IND_{\{R_j(f - \pi f,\kappa_0 r/2) \geq m\rho_0/4\}} \leq \frac{\tau}{12}$.
\end{description}
On this event, for every $f \in \F_1$ there are at least $(1-\tau/4)n$  blocks $I_j$ with the following properties:
\begin{description}
\item{$\bullet$} $|\pi f (X_i)-f^*(X_i)| \geq \kappa_0 r $ on at least $m\rho_0/2$ coordinates in $I_j$, and
\item{$\bullet$} $|f (X_i) - \pi f (X_i)| \geq \kappa_0 r/2$ on at most $m\rho_0/4$ coordinates in $I_j$.
\end{description}
Hence, in each one of the $(1-\tau/4)n$ well-behaved blocks there are at least $m \rho_0/4$ coordinates $X_i$ that satisfy
$$
|f (X_i)-f^*(X_i)| \geq |\pi f(X_i) - f^*(X_i)| - |f (X_i)-\pi f(X_i)| \geq \kappa_0 r/2~.
$$
In particular, since $\|f-f^*\|_{L_2}=r$, one has
\begin{equation} \label{eq:quad-lower-in-proof}
{\Q}_{f,f^*}(j) \geq (\rho_0 \kappa_0^2/16)\|f-f^*\|_{L_2}^2~.
\end{equation}
Moreover, \eqref{eq:quad-lower-in-proof} is positive homogeneous in $f-f^*$ and $\F_1$ is star-shaped around $f^*$, implying that \eqref{eq:quad-lower-in-proof} holds for every $f \in \F_1$ as long as $\|f-f^*\|_{L_2} \geq r$.

All that remains is to establish that $(B)$ holds with a sufficiently high probability. To that end, let
$$
\Psi(X_1,...,X_N)=\sup_{f \in \F_1} \frac{1}{n}\sum_{i =1}^n \IND_{\{R_j(f-\pi f,\kappa_0 r/2) \geq m\rho_0/4\}}~.
$$
First note that, by the bounded differences inequality (see, e.g., \cite{BoLuMa13}), with probability at least $1-\exp(-c_3u^2)$,
$$
\Psi(X_1,...,X_N) \leq \E \Psi + \frac{u}{\sqrt{n}}~.
$$
Therefore, if $\E \Psi \leq \tau/24$ and $u \leq \sqrt{n} \tau/24$, then with probability at least $1-\exp(-c_4\tau^2n)$, $\Psi(X_1,...,X_N) \leq \tau/12$ as required.

Thus, the final step in the proof is to show that
$$
(*)=\E \sup_{f \in \F_1} \frac{1}{n}\sum_{i =1}^n \IND_{\{R_j(f-\pi f,\kappa_0 r/2) \geq m\rho_0/4\}} \leq \frac{\tau}{24}~.
$$

Using that, for $\alpha>0$, $\IND_{\{|x| \geq \alpha\}} \leq \alpha^{-1} |x|$,
$$
\IND_{\{R_j(f-\pi f,\kappa_0 r/2) \geq m\rho_0/4\}} \leq \frac{4}{\rho_0 m} R_j(f-\pi f,\kappa_0 r/2) = \frac{4}{\rho_0 m} \sum_{i \in I_j} \IND_{\{|f (X_i)-\pi f(X_i)| \geq \kappa_0 r/2\}}~.
$$
%Recall that $(X_i)_{i=1}^N$ are i.i.d.\ random variables and that $nm=N$.
%By the Gin\'{e}-Zinn symmetrization Theorem \cite{ }, followed by the
%contraction principle for Bernoulli processes \cite{ }, applied to
%the Lipschitz function $\phi(t)=|t|$,
Thus, recalling that $nm=N$,
\begin{align*}
(*) \leq & \frac{4}{\rho_0} \cdot  \E \sup_{f \in \F_1}
\frac{1}{N}\sum_{i=1}^N \IND_{\{|f-\pi f|(X_i) \geq \kappa_0 r/2\}}
\\
 \leq & \frac{8}{\rho_0 \kappa_0 r} \E \sup_{f \in \F_1}
 \frac{1}{N}\sum_{i=1}^N |f (X_i)-\pi f(X_i)|
\quad \text{(using $\IND_{\{|x| \geq \alpha\}} \leq \alpha^{-1} |x|$ again)}
\\
\leq & \frac{8}{\rho_0 \kappa_0 r} \left( \E \sup_{f \in \F_1}
  \frac{1}{N}\sum_{i=1}^N \left(|f (X_i)-\pi f(X_i)| - \E |f(X)-\pi
    f(X)|\right)
+ \sup_{f \in \F_1} \E|f(X)-\pi f(X)| \right)
\\
\leq & \frac{16}{\rho_0 \kappa_0 r}  \E \sup_{f \in \F_1}
\left|\frac{1}{N}\sum_{i=1}^N \eps_i (f(X_i)-\pi f(X_i)) \right| +
\frac{8}{\rho_0 \kappa_0 r} \cdot \eta r \\
& \text{(where $(\eps_i)_{i=1}^N$ are independent symmetric random signs)}~.
\end{align*}
In the last step we used a standard symmetrization argument, see,
for example, \cite{deGi99}.
To conclude the proof, observe that $({8}/{\rho_0 \kappa_0 r}) \cdot \eta r \leq \tau/48$ when
\begin{equation} \label{eq:eta-cond-1}
\eta \leq c_5 \rho_0 \kappa_0 \tau
\end{equation}
and for a suitable absolute constant $c_5$. The fact that
\begin{equation} \label{eq:cond-01}
\frac{16}{\rho_0 \kappa_0 r}  \E \sup_{f \in \F_1} \left|\frac{1}{N}\sum_{i=1}^N \eps_i (f-f^*)(X_i) \right| \leq \frac{\tau}{48}
\end{equation}
follows because $r \geq 2r^*$ and invoking \eqref{eq:r*-emp}, as long as
\[
%begin{equation} \label{eq:kappa1-cond1}
\kappa_1 \leq c_6 \rho_0 \kappa_0 \tau
\]
for an absolute constant $c_6$.
\endproof

\subsubsection*{The multiplier component} \label{sec:multi}
Next we complete the proof of Lemma \ref{lemma:components-of-main} by showing that, with high probability, if $\|f-f^*\|_{L_2} \geq r$, then
$$
{\M}_{f,f^*}(j) \leq -(3C_1/4) \|f-f^*\|_{L_2}^2
$$
on at most $\tau n$ blocks. In the first step towards this goal, we
consider a single function:
\begin{Lemma} \label{lemma:BE}
There exists an absolute constant $C_2$ for which the following
holds. Assume  that for every $f \in \F$, $\E (f^*(X)-Y)(f(X)-f^*(X)) \geq 0$.
 If $f \in \F_1$, then with probability at least $1-2\exp(-C_2\tau^2n)$,
$$
\left|\left\{j : {\M}_{f,f^*}(j) \geq -\frac{C_1r^2}{2}\right\} \right| \geq (1-\tau/8)n~.
$$
In particular, the same assertion holds uniformly for any fixed subset
$H_2 \subset \F_1$ of cardinality at most $\exp(C_2\tau^2n/2)$.
\end{Lemma}

Before proving Lemma \ref{lemma:BE}, consider the assumption that
\begin{equation} \label{eq:conv-cond}
\E (f^*(X)-Y)(f(X)-f^*(X)) \geq 0 \ \ \ {\rm for \ every \ } \ \ f \in \F~.
\end{equation}
It is straightforward to verify that \eqref{eq:conv-cond} is satisfied
under the assumptions of Theorems \ref{thm:main-general} and
\ref{thm:main-independent}.  Indeed, if $\F$ is closed and convex then
\eqref{eq:conv-cond} is just the characterization of $f^*$ as the
nearest point to $Y$ in $\F$ in the $L_2$ sense. On the other hand, if $Y=f_0(X)+W$
for $f_0 \in \F$ and $W$ that is mean-zero and independent of $X$,
then $\E (f^*(X)-Y)(f(X)-f^*(X))=-\E W(f(X)-f_0(X)) =0$ for every $f \in
\F$. In fact, \eqref{eq:conv-cond} is the only structural assumption
on the `location' of $Y$ relative to $\F$ that is required for our
analysis.

Another observation is that \eqref{eq:conv-cond} passes to ${\rm
  star}(\F,f^*)$, simply because any function in that set is of the
form $h=\lambda f + (1-\lambda)f^*$ for some $f \in \F$ and $0 \leq
\lambda \leq 1$ and therefore, $(f^*(X)-Y)(h(X)-f^*(X))=\lambda (f^*(X)-Y)(f(X)-f^*(X))$.

\medskip

\noindent{\bf Proof of Lemma \ref{lemma:BE}.}
Set $\xi=f^*(X)-Y$, put $U=\xi (f(X)-f^*(X))$ and observe that by \eqref{eq:conv-cond}, $\E U \geq 0$. Also note that for every $j$,
$$
{\cal M}_{f,f^*}(j)=\frac{1}{m}\sum_{i \in I_j} \xi_i (f(X_i)-f^*(X_i))=\frac{1}{m} \sum_{i \in I_j} U_i~,
$$
for $(U_i)_{i=1}^M$ that are independent copies $U$. It is
straightforward to verify that,
with $(\eps_i)_{i=1}^N$ defined as independent symmetric random signs,
\begin{align*}
\PROB \left( \left|\frac{1}{m}\sum_{i=1}^m U_i - \E U \right|  \geq  t
\right) \leq &
\frac{\E \left|\sum_{i=1}^m (U_i - \E U) \right|}{mt} \leq 2\frac{\E \left|\sum_{i=1}^m \eps_i U_i\right|}{mt}
\\
\leq & 2 \frac{\|U\|_{L_2}}{\sqrt{m} t} = 2\sqrt{n} \frac{\|U\|_{L_2}}{\sqrt{N} t} = (*)~.
\end{align*}
By the norm equivalence assumption of Theorem \ref{thm:main-general},
$$
\|U\|_{L_2} \leq \|\xi\|_{L_4} \|f-f^*\|_{L_4} \leq L^2 \sigma r~,
$$
and by the independence assumption of Theorem \ref{thm:main-independent},
$$
\|U\|_{L_2} = \|\xi\|_{L_2} \|f-f^*\|_{L_2} = \sigma r~.
$$
Hence, setting $t=(C_1/2) r^2$ and noting that
$$
\sqrt{\frac{n}{N}} \leq \sqrt{\theta} \frac{r}{\sigma}~,
$$
it follows that
$$
(*) \leq 2L^2  \frac{r \sigma\sqrt{n}}{t \sqrt{N}} = \frac{4L^2}{C_1} \frac{\sqrt{n} \sigma}{\sqrt{N}r} \leq \frac{\tau}{16}
$$
whenever
\[
%begin{equation} \label{eq:theta-cond-1}
\theta \leq (C_1 \tau /64L^2)^2~.
\]
Therefore, with probability at least $1-\tau/16$,
$$
\frac{1}{m} \sum_{i=1}^m U_i \geq \E U - C_1 r^2/2 \geq - C_1 r^2/2~.
$$

Finally, consider the independent Bernoulli random variables
$(\IND_{\{\M_{f,f^*}(j) \leq -C_1 r^2/2\}})_{j=1}^n$ that have mean at
most $\tau/16$. By concentration of Binomial random variables,
there is an absolute
constant $C_2$ such that, with probability at least
$1-2\exp(-C_2\tau^2 n)$, there are at least $(1-\tau/8)$
blocks that satisfy $\M_{f,f^*}(j) \geq -C_1r^2/2$, as claimed.
\endproof

Just as before, one may select any fixed subset $H_2 \subset \F_1$ of
cardinality $\exp(C_2\tau^2n/2)$, and the assertion of Lemma
\ref{lemma:BE} holds with high probability and uniformly for every $h
\in H_2$. The choice of $H_2$ requires some care. It cannot be just an
arbitrary maximal separated set.

\begin{Lemma} \label{lemma:the-net}
There exists a subset $H_2 \subset \F_1$ of cardinality at most
$\exp(C_2n\tau^2/2)$ such that for every $f \in \F_1$ there is some $h \in H_2$ that satisfies
$$
\|f-h\|_{L_2} \leq 2 \eta r, \ \ {\rm and} \ \ \E \xi(f(X)-h(X)) \geq 0~.
$$
\end{Lemma}

\proof Let $H^\prime$ be an maximal $\eta r$-separated subset of $\F_1$. Recall that by the choice of $r$ and \eqref{eq:r*-entropy}, $\log |H^\prime| \leq  \kappa_3^2 N\min\{1,\sigma^{-2} r^2\}$, which is smaller than $(C_2/2)\tau^2 n$ when
\[
%begin{equation} \label{eq:kappa3-cond-2}
\kappa_3^2 \leq (C_2/2) \theta \tau^2 .
\]

The class $\F$ is a locally compact subset of $L_2(\mu)$ and therefore, so is
$\F_1$.
Thus, for every $h^\prime \in H^\prime$, the intersections of $\F_1$ with the $L_2$ balls $B(h^\prime,\eta r)$ are compact and the continuous linear functional on $L_2$
$$
f \to \E \xi f(X)
$$
attains its minimum in each one of the sets $\F_1 \cap B(h^\prime,\eta r)$. Let $h$ be such a minimizer and set $H_2$ to be the collection of these minimizers. Hence, for every $f \in \F_1$ there is some $h \in H_2$ for which $\|f-h\|_{L_2} \leq 2 \eta r$ and $\E f(X)\xi \geq \E h(X)\xi$, as claimed.
\endproof

With Lemma \ref{lemma:the-net} at our disposal, for every $f \in \F_1$
define $\pi f \in H_2$ for which $\E \xi(f(X)-\pi f(X)) \geq 0$
and $\|f-\pi f\|_{L_2} \leq 2\eta r$, as above.
By Lemma \ref{lemma:BE}, with probability at least $1-2\exp(-C_2\tau^2
n/2)$, for every $h \in H_2$ there are at least $(1-\tau/8)n $
blocks $I_j$ with
$$
{\M}_{h,f^*}(j) = \frac{2}{m} \sum_{i \in I_j} \xi_i (h(X_i)-f^*(X_i)) \geq -C_1r^2/2~.
$$
Hence, to complete the proof it suffices to show that for every $f \in
\F_1$, there are at most $(7/8)n\tau$ blocks $I_j$ with
$$
{\M}_{f,f^*}(j) \leq {\M}_{\pi f,f^*}(j) - C_1r^2/4~.
$$
Indeed, on that event, for every $f \in \F_1$ there are at most $n\tau$ blocks $I_j$ that satisfy ${\M}_{f,f^*}(j) \leq -(3/4)C_1r^2$, as required.

Again, establishing this estimate for $\F_1$ yields that on the same
event, if $\|f-f^*\|_{L_2}^2 \geq r^2$, then ${\M}_{f,f^*}(j) \leq
-(3/4)C_1\|f-f^*\|_{L_2}^2$ on at most $n\tau$ coordinate blocks,
since $\M$ is homogeneous in
 $f-f^*$ and $\F_1$ is star-shaped around $f^*$.

\begin{Lemma} \label{lemma:oscillations}
There exists an absolute constant $c$ for which, with probability at least $1-2\exp(-c\tau^2n)$,
$$
\sup_{f \in \F_1} \frac{1}{n}\sum_{j=1}^n \IND_{\{{\M}_{f,f^*}(j) - {\M}_{\pi f,f^*}(j) \leq -C_1r^2/4\}} \leq \frac{7\tau}{8}~.
$$
\end{Lemma}

\proof Recall that by the definition of $\pi f$, $\E \xi f(X) \geq \E \xi \pi f(X)$. Therefore,
$$
\E ({\M}_{f,f^*}-{\M}_{\pi f, f^*})= 2\E\xi(f(X)-\pi f(X)) \geq 0~.
$$
To simplify notation, set
$$
W_{f,\pi f}(j) = {\M}_{f,f^*}(j)-{\M}_{\pi f, f^*}(j)= \frac{1}{m}\sum_{i \in I_j} \xi_i (f(X_i)-\pi f(X_i))~,
$$
and note that $\E W_{f,\pi f}(j) \geq 0$.

Consider the supremum of Binomial random variables
$$
\Psi = \sup_{f \in \F_1} \frac{1}{n}\sum_{j=1}^n \IND_{\{W_{f,\pi f}(j) \leq -C_1r^2/4\}}~.
$$
By the bounded differences inequality, with probability at least
$1-2\exp(-cu^2)$,
$\Psi \leq \E \Psi + u/\sqrt{n}$. Setting $u=(6/8)\sqrt{n}\tau$, all that remains is to show that
$$
\E \Psi = \E \sup_{f \in \F_1} \frac{1}{n}\sum_{j=1}^n \IND_{\{W_{f,\pi f}(j) \leq -C_1r^2/4\}} \leq \frac{\tau}{8}~.
$$
Since $\E W_{f,\pi f} \geq 0$, using $\IND_{\{|x| \geq \alpha\}} \leq \alpha^{-1} |x|$,
\begin{align*}
 \E \sup_{f \in \F_1} \frac{1}{n}\sum_{j=1}^n \IND_{\{W_{f,\pi f}(j)
   \leq -C_1r^2/4\}} & \leq \E \sup_{f \in \F_1} \frac{1}{n}\sum_{j=1}^n \IND_{\{W_{f,\pi f}(j)-\E W_{f,\pi f}(j)  \leq -C_1r^2/4\}}
\\
& \leq  \E \sup_{f \in \F_1} \frac{1}{n}\sum_{j=1}^n \IND_{\{|W_{f,\pi
    f}(j)-\E W_{f,\pi f}(j)|  \geq C_1r^2/4\}} \\
& \leq \frac{4}{C_1r^2} \E \sup_{f \in \F_1} \frac{1}{n} \sum_{j=1}^n |W_{f,\pi f}(j)-\E W_{f,\pi f}(j)|~.
\end{align*}
The next step is to centre the process. We show that for every $f \in \F_1$, the centring term satisfies
$$
\frac{4}{C_1r^2} \E |W_{f,\pi f}(j)-\E W_{f,\pi f}(j)| \leq \frac{\tau}{16}~.
$$
Indeed, symmetrizing, applying either the assumption of norm
equivalence of Theorem \ref{thm:main-general} or the independence
assumption of Theorem \ref{thm:main-independent},
 and recalling that $\|f-\pi f\|_{L_2} \leq 2\eta r$, it is evident that
\begin{align*}
\frac{4}{C_1r^2}\E |W_{f,\pi f} - \E W_{f,\pi f} | \leq &
\frac{8}{C_1r^2} \E \left|\frac{1}{m}\sum_{i=1}^m \eps_i \xi_i (f(X_i)
  - \pi f(X_i)) \right|
\\
\leq & \frac{16L^2}{C_1} \frac{\sigma \eta r}{r^2 \sqrt{m}}
\\
=&  \frac{16L^2}{C_1} \frac{\eta \sigma  \sqrt{n} }{r \sqrt{N}} \leq \frac{16L^2}{C_1} \eta \sqrt{\theta}~,
\end{align*}
again using the fact that $\sqrt{n/N} \leq \sqrt{\theta}(r/\sigma)$. Clearly,
$$
\frac{16L^2}{C_1} \eta \sqrt{\theta} \leq \frac{\tau}{16}
$$
when
\[
%begin{equation} \label{eq:theta-cond-2}
\theta \leq (C_3 \tau /L^2 \eta )^2
\]
for a suitable absolute constant $C_3$.

Finally, we need to bound the centred empirical process. This is done by
standard techniques of symmetrization, the contraction theorem for
Bernoulli processes (for the Lipschitz function $\phi(t)=|t|$)
  and then de-symmetrization (see, e.g., \cite{LeTa91}):
\begin{eqnarray*}
\lefteqn{
 \E \sup_{f \in \F_1} \left| \frac{1}{n} \sum_{j=1}^n  \left|W_{f,\pi
       f}(j) - \E W_{f,\pi f} \right| - \E \left|W_{f,\pi f}(j) - \E
     W_{f,\pi f} \right| \right|  } \\
& \leq & 2\E \sup_{f \in \F_1} \left| \frac{1}{n} \sum_{j=1}^n
  \eps_j\left|W_{f,\pi f}(j) - \E W_{f,\pi f} \right| \right|
\\
& \leq &  2 \E \sup_{f \in \F_1} \left| \frac{1}{n} \sum_{j=1}^n  \eps_j\left(W_{f,\pi f}(j) - \E W_{f,\pi f} \right) \right|
\\
& \leq & 4 \E \sup_{f \in \F_1} \left| \frac{1}{n} \sum_{j=1}^n  \left(W_{f,\pi f}(j) - \E W_{f,\pi f} \right)\right|~.
\end{eqnarray*}
Since
$$
W_{f,\pi f}(j) = \frac{1}{m} \sum_{i \in I_j} \xi_i (f(X_i)-\pi f(X_i))
$$
one has that
\begin{align*}
\frac{16}{C_1r^2} & \E \sup_{f \in \F_1} \left|
  \frac{1}{N}\sum_{i=1}^N \xi_i (f(X_i)-\pi f(X_i)) - \E \xi (f(X)-\pi
  f(X)) \right| \nonumber
\\
\leq & \frac{32}{C_1r^2} \E \sup_{f \in \F_1} \left| \sum_{i=1}^N \eps_i \xi_i (f(X_i)-f^*(X_i)) \right| \leq \frac{\tau}{16}
\end{align*}
provided that
\begin{equation} \label{eq:cond-3}
\E \sup_{f \in \F_1} \left| \frac{1}{\sqrt{N}} \sum_{i=1}^N \eps_i \xi_i (f(X_i)-f^*(X_i)) \right| \leq C_4 \tau \sqrt{N} r^2~,
\end{equation}
that is, when
\[
%begin{equation} \label{eq:kappa2-cond1}
\kappa_2 \leq C_4 \tau~,
\]
concluding the proof.
\endproof

\subsection{Champions league--proof} \label{sec:champions-proof}

Finally, it remains to prove Proposition \ref{thm:winners}.  Recall
that, with probability at least $1-2\exp\bigl(-c_0N
\min\{1,\sigma^{-2}r^2\}\bigr)$ with
respect to the sample $(X_i,Y_i)_{i=1}^{2N}$, we have been able to
identify a set of ``qualifiers" that have not lost a single match
in the preliminary round subset; that is, $H \subset \F$, consisting of $f^*$ and possibly
other functions that satisfy $\|f-f^*\|_{L_2} \leq (\beta/\alpha) r$.
In the rest of this section we work conditionally on this ``good''
event.

While producing a function that is close to
$f^*$ solves the estimation problem, the question of prediction
requires an additional step: we would like to choose one of the
qualifiers that has an almost optimal statistical performance: a function $\wh{f}$ that satisfies
\begin{equation} \label{eq:oracle-champions}
\E \bigl( (\wh{f}(X)-Y)^2 | (X_i,Y_i)_{i=2N+1}^{3N}\bigr) \leq \E(f^*(X)-Y)^2 + Cr^2
\end{equation}
for an appropriate $C$. We show that this is possible with the
required high probability for $C=16(\beta/\alpha)^2$.
To this end, set $r_1=(\beta/\alpha)r$.

Recall that for $f,h \in {\rm star}(\F,f^*)$, $\Psi_{h,f} = (h(X)-f(X))(f(X)-Y)$ and note that
$$
{\M}_{f,f^*} = \frac{2}{m}\sum_{i=1}^m \Psi_{f,f^*}(X_i,Y_i)~.
$$
Also, as noted previously, since $\E \Psi_{f,f^*} \geq 0$ for every $f \in \F$, the same holds for any $f \in {\rm star}(\F,f^*)$.

As was stated in Section \ref{sec:winners}, the sub-sample
$(X_i,Y_i)_{i=2N+1}^{3N}$ is used to identify a function $\wh{f} \in \F$ for
which $\E \Psi_{f^*,\wh{f}} \geq -2r_1^2$. By Lemma
\ref{lemma:good-f-oracle}, this indeed suffices to establish
\eqref{eq:oracle-champions} for $C=16(\beta/\alpha)^2$.

The ``champions league" round is designed to have ``home-and-away" legs. For the partition $(I_j)_{j=1}^n$ of $(X_i,Y_i)_{i=2N+1}^{3N}$, $f$ wins its home match against $h$ if
$$
\frac{2}{m} \sum_{i \in I_j} \Psi_{h,f}(X_i,Y_i) \geq - r^2_1/10
$$
for more than half of the blocks $I_j$. We show that
$f^*$ wins all of its home matches, implying that the set of possible
choices of $\wh{f}$ is nonempty, and that if $\E \Psi_{f^*,f} \leq
-2r^2_1$, then $f$ loses its home match against $f^*$. On that event,
a function that wins all of its home matches must satisfy that
$\E\Psi_{f^*,\wh{f}} \geq -2r_1^2$, and this observation concludes the
proof of Theorem \ref{thm:main-general}.

The main ingredient in the proof is the following lemma.
\begin{Lemma} \label{lemma:multi-for-champ}
Let $\F_2={\rm star}(\F,f^*) \cap (f^*+r_1 D)$. Under the conditions of
Theorem \ref{thm:main-general}, there is an absolute constant $c$, such that, with probability at least $1-2\exp(-cn)$,
$$
\sup_{f \in \F_2} \frac{1}{n} \sum_{j=1}^n \IND_{\{|{\M}_{f,f^*}(j) - \E {\M}_{f,f^*} | \geq r_1^2/10 \} } < \frac{1}{2}~.
$$
\end{Lemma}

\proof Using the assumption of norm equivalence of Theorem
\ref{thm:main-general}, one may verify that
$\|\Psi_{f,f^*}\|_{L_2} \leq L^2 r_1 \sigma$. Therefore,
\begin{equation} \label{eq:mean-in-proof-champs}
\E \left|{\M}_{f,f^*} - \E {\M}_{f,f^*}\right| \leq \frac{2}{\sqrt{m}} \|\Psi_{f,f^*}\|_{L_2} \leq \frac{2}{\sqrt{m}} \cdot L^2 r_1\sigma~.
\end{equation}
Setting
$$
Z=\sup_{f \in \F_2} \frac{1}{n} \sum_{j=1}^n \IND_{\{|{\M}_{f,f^*}(j) - \E {\M}_{f,f^*} | \geq r_1^2/10 \} },
$$
it follows that $\E Z$ is at most
\begin{align*}
 \E \sup_{f \in \F_2} \frac{10}{n r_1^2} \sum_{j=1}^n & \left(| {\M}_{f,f^*}(j) - \E {\M}_{f,f^*}| - \E | {\M}_{f,f^*}(j) - \E {\M}_{f,f^*}|\right)
\\
+ & \frac{10}{n r_1^2} \cdot \sup_{f \in \F_2} \E | {\M}_{f,f^*}(j) - \E {\M}_{f,f^*}| = (I)+(II)~.
\end{align*}
Applying the same argument used in the previous section--namely, symmetrization, followed by contraction for a Bernoulli process and the Lipschitz function $\phi(t)=|t|$ and de-symmetrization--one has that for absolute constants $c_1,c_2$, and $c_3$,
\begin{equation*}
(I) \leq \frac{c_1}{r_1^2} \cdot \E  \sup_{f \in \F_2}  \left|\frac{1}{n} \sum_{j=1}^n {\M}_{f,f^*}(j) - \E {\M}_{f,f^*}(j) \right|
\leq \frac{c_2}{r_1^2} \cdot \E  \sup_{f \in \F_2}  \left|\frac{1}{N} \sum_{i=1}^N \eps_i \xi_i (f-f^*)(X_i) \right| \leq \frac{1}{8}
\end{equation*}
when
\[
%begin{equation} \label{eq:kappa2-cond-2}
\kappa_2 \leq c_3~.
\]
Also, by \eqref{eq:mean-in-proof-champs},
$$
(II) \leq \frac{10 L^2 \sigma}{r_1} \cdot \sqrt{\frac{n}{N}} \leq \frac{1}{8}~,
$$
provided that $n \leq c_4 L^4 \cdot (r_1^2/\sigma^2) N$ for an
absolute constant $c_4$.
Since  $r_1 =2(\beta/\alpha)r$, it suffices that
\[
\theta \leq c_4 L^4~.
\]
Thus, $\E Z \leq 1/4$, and by the bounded differences inequality applied to $Z$, one has that with probability at least $1-2\exp(-cu^2)$, $Z \leq \E Z +u/\sqrt{n}$. The claim follows by selecting $u=\sqrt{n}/8$.
\endproof

Recall that for every $f \in H$, $\|f-f^*\|_{L_2} \leq
(\beta/\alpha)r$ and consider the ``good'' event from Lemma
\ref{lemma:multi-for-champ}.
For any $f \in \F_2$, and, in particular, for any qualifier $f \in H$,
\begin{equation} \label{eq:champions-good-event}
\E {\M}_{f,f^*} - r_1^2/10 \leq {\M}_{f,f^*}(j) \leq \E {\M}_{f,f^*} + r_1^2/10
\end{equation}
on more than $n/2$ blocks. Moreover, if $\|f-f^*\|_{L_2} \leq r_1$ and $\E \Psi_{f^*,f} \leq -2 r_1^2$, then $\E \Psi_{f,f^*} \geq r_1^2$. Indeed,
\begin{align*}
-2r_1^2 \geq & \E (f^*(X)-f (X)) \cdot (f(X)-Y)
\\
= & -\|f-f^*\|_{L_2}^2 +\E (f^*(X)-f (X)) \cdot (f^*(X)-Y)
\\
\geq & -r^2 +\E (f^*(X)-f (X)) \cdot (f^*(X)-Y) = -r_1^2 - \E\Psi_{f,f^*}~.
\end{align*}
Therefore, $\E \M_{f,f^*} \geq 2r_1^2$, and on the event \eqref{eq:champions-good-event},
$$
\M_{f,f^*}(j) \geq \E \M_{f,f^*} - r_1^2/10 \geq r_1^2~.
$$
Finally, since $\Psi_{f^*,f} = -(f^*(X)-f(X))^2 - \Psi_{f,f^*}$, on that event and the same coordinate blocks,
$$
\frac{1}{m} \sum_{i \in I_j} \Psi_{f^*,f}(j) = -\frac{2}{m} \sum_{i=1}^m (f^*(X_i)-f(X_i))^2 - {\M}_{f,f^*} \leq -r_1^2~.
$$
Thus, $f$ is defeated by $f^*$ in the majority of the blocks, and in
particular, loses its home match against $f^*$.  It follows that, with
probability at least $1-2\exp(-cn)$, any function $\wh{f}$ selected in
the champions league round must satisfy that $\Psi_{f^*,\wh{f}} \geq
-2r_1^2$, and by Lemma \ref{lemma:good-f-oracle},
$$
\E\bigl((\wh{f}(X)-Y)^2 | (X_i,Y_i)_{2N+1}^{3N} \bigr) \leq \E(f^*(X)-Y) + 4r_1^2~.
$$
\endproof

\section*{Additional remarks}

%SHAHAR: I THINK THAT THESE SHOULD BE MOVED TO A ``CONCLUDING REMARKS" SECTION -- IF AT ALL.
%
It should be noted that the difference between $r^*(f^*)$ and $r^*$ is not a major
issue in most interesting cases. $\lambda_{\Q}(\kappa,\eta,h)$, $\lambda_{\M}(\kappa,\eta,h)$, and $r_E(\kappa,h)$
typically do not depend heavily on $h$ and little is lost by taking
the supremum over $h\in \F$. $\wt{r}_{\M}(c_1,\sigma)$ becomes a better estimate of $\ol{r}_{\M}(c_1,f^*)$
as $\sigma$ gets closer to $\|f^*-Y\|_{L_2}$. Even if this value is not
known in advance, it is easy to design a two-stage procedure that constructs
a data-dependent estimate of $\|f^*-Y\|_{L_2}$ and then uses the procedure
of Theorem \ref{thm:main-general} with a tight upper bound for the value of $\sigma$ obtained from the first
 stage. In particular, it is easy to find a value of $\sigma$ that is
 within a constant factor of the optimum.
This is an issue of secondary importance and we omit the straightforward details.

\vskip0.4cm

The major problem that remains open is the identity of the accuracy edge. To date, there is a single generic example in which one may attain an accuracy smaller than  $Cr^*$ and in that case the accuracy attained is proportional to $\lambda^*$ and with the optimal confidence at that level, namely,
$1-2\exp \bigl(-c N\min\{1,\sigma^{-2} \lambda_{\M}^2(\kappa_2/\sigma,\eta_2)\}\bigr)$.

This fact has recently been established in \cite{Men16} in a very special situation: when $\F$ is a convex, $L$-sub-Gaussian class of functions; \emph{all} the admissible targets are of the form
$Y=f_0(X)+W$ for $f_0 \in \F$ and $W$ that is sub-Gaussian, zero-mean, independent of $X$. The procedure used is a modification of ERM: one replaces $\F$ with an appropriate net, thus `erasing' all the fine structure of $\F$ at the right level, and then runs ERM on the net. The idea behind this procedure is straightforward: if one is interested in accuracy $r$, one is insensitive to perturbations of that order. From that perspective, a net with a mesh that is proportional to $r$ is just as good as the entire class. It turns out that $r^*$ of the net is proportional to $\lambda^*$ the original class $\F$.
Unfortunately, all the highly restrictive assumptions are essential to the proof and cannot be relaxed at all.

It should be noted that the median-of-means tournament may be modified
in exactly the same way as ERM is modified in \cite{Men16}, leading to
an accuracy that is proportional to $\lambda^*$ when $\F$ is a convex,
$L$-sub-Gaussian class and for an independent noise $W$ that may be
heavy-tailed. We decided not pursue this point further because it is a
very special case, and shifts the emphasis of the article from the
question of the optimal accuracy/confidence tradeoff to the nature of
the accuracy edge, a problem of a different nature.

As far as the latter is concerned, it is not clear whether the gap
between $\lambda^*$ and $r^*$ can be closed in other cases. It is
possible that both independent noise and a sub-Gaussian class are
essential in attaining an accuracy proportional to $\lambda^*$, and
under weaker assumptions, the true accuracy edge lies somewhere between
$\lambda^*$ and $r^*$. We leave that question for future study.

Finally, it should be noted that in this work we completely ignore the
algorithmic aspects of the new procedure. While computing the
empirical risk minimizer often leads to thoroughly studied and well
understood convex
optimization problems, finding the winner of the median-of-means
tournament in a computationally efficient manner is a highly
nontrivial--and perhaps not hopeless--problem that goes beyong the
scope of this paper. Techniques of ``derivative-free optimization''
with ``function comparison oracle'' may be useful, see, for example,  \cite{JaNoRe12}.

\bibliographystyle{plain}
%\bibliography{tournament}

\begin{thebibliography}{10}

\bibitem{AMS02}
N.~Alon, Y.~Matias, and M.~Szegedy.
\newblock The space complexity of approximating the frequency moments.
\newblock {\em Journal of Computer and System Sciences}, 58:137--147, 2002.

\bibitem{AnBa99}
M.~Anthony and P.~L. Bartlett.
\newblock {\em Neural network learning: theoretical foundations}.
\newblock Cambridge University Press, Cambridge, 1999.

\bibitem{AuCa11}
J.-Y. Audibert and O.~Catoni.
\newblock Robust linear least squares regression.
\newblock {\em The Annals of Statistics}, 39:2766--2794, 2011.

\bibitem{BGMN05}
F.~Barthe, O.~Gu{\'e}don, S.~Mendelson, and A.~Naor.
\newblock A probabilistic approach to the geometry of the {$l^n_p$}-ball.
\newblock {\em Ann. Probab.}, 33(2):480--513, 2005.

\bibitem{BaKo03}
F.~Barthe and A.~Koldobsky.
\newblock Extremal slabs in the cube and the {L}aplace transform.
\newblock {\em Adv. Math.}, 174(1):89--114, 2003.

\bibitem{BaBoMe05}
P.L. Bartlett, O.~Bousquet, and S.~Mendelson.
\newblock Localized {R}ademacher complexities.
\newblock {\em Annals of Statistics}, 33:1497--1537, 2005.

\bibitem{BoLuMa13}
S.~Boucheron, G.~Lugosi, and P.~Massart.
\newblock {\em Concentration inequalities: A Nonasymptotic Theory of
  Independence}.
\newblock Oxford University Press, 2013.

\bibitem{BrJoLu15}
C.~Brownlees, E.~Joly, and G.~Lugosi.
\newblock Empirical risk minimization for heavy-tailed losses.
\newblock {\em Annals of Statistics}, 43:2507--2536, 2015.

\bibitem{BuCeLu13}
S.~Bubeck, N.~Cesa-Bianchi, and G.~Lugosi.
\newblock Badits with heavy tail.
\newblock {\em IEEE Transactions on Information Theory}, 59:7711--7717, 2013.

\bibitem{BuvdG11}
P.~B{\"u}hlmann and S.~{van de Geer}.
\newblock {\em Statistics for high-dimensional data}.
\newblock Springer Series in Statistics. Springer, Heidelberg, 2011.
\newblock Methods, theory and applications.

\bibitem{Cat10}
O.~Catoni.
\newblock Challenging the empirical mean and empirical variance: a deviation
  study.
\newblock {\em Annales de l'Institut Henri Poincar{\'e}, Probabilit{\'e}s et
  Statistiques}, 48(4):1148--1185, 2012.

\bibitem{deGi99}
V.H. de~la Pe\~na and E.~Gin\'e.
\newblock {\em Decoupling: from Dependence to Independence}.
\newblock Springer, New York, 1999.

\bibitem{DGL96}
L.~Devroye, L.~Gy{\"o}rfi, and G.~Lugosi.
\newblock {\em A probabilistic theory of pattern recognition}, volume~31 of
  {\em Applications of Mathematics (New York)}.
\newblock Springer-Verlag, New York, 1996.

\bibitem{DeLeLuOl16}
L.~Devroye, M.~Lerasle, G.~Lugosi, and R.I. Oliveira.
\newblock Sub-gausssian mean estimators.
\newblock {\em Annals of Statistics}, 2016.

\bibitem{Dud:book}
R.~M. Dudley.
\newblock {\em Uniform central limit theorems}, volume 142 of {\em Cambridge
  Studies in Advanced Mathematics}.
\newblock Cambridge University Press, New York, second edition, 2014.

\bibitem{FoRa13}
S.~Foucart and H.~Rauhut.
\newblock {\em A mathematical introduction to compressive sensing}.
\newblock Applied and Numerical Harmonic Analysis. Birkh\"auser/Springer, New
  York, 2013.

\bibitem{GyKoKrWa02}
L.~Gy\"orfi, M.~Kohler, A.~Krzy\.zak, and H~Walk.
\newblock {\em A Distribution-Free Theory of Nonparametric Regression}.
\newblock Springer, New York, 2002.

\bibitem{HsuSa13}
D.~Hsu and S.~Sabato.
\newblock Approximate loss minimization with heavy tails.
\newblock {\em Computing Research Repository}, abs/1307.1827, 2013.

\bibitem{HsSa16}
D.~Hsu and S.~Sabato.
\newblock Loss minimization and parameter estimation with heavy tails.
\newblock {\em Journal of Machine Learning Research}, 17:1--40, 2016.

\bibitem{JaNoRe12}
K.~Jamieson, R.~Nowak, and B.~Recht.
\newblock Query complexity of derivative-free optimization.
\newblock In {\em Advances in Neural Information Processing Systems}, pages
  2672--2680, 2012.

\bibitem{JeVaVa86}
M.~Jerrum, L.~Valiant, and V.~Vazirani.
\newblock Random generation of combinatorial structures from a uniform
  distribution.
\newblock {\em Theoretical Computer Science}, 43:186--188, 1986.

\bibitem{Kolt08}
V.~Koltchinskii.
\newblock {\em Oracle inequalities in empirical risk minimization and sparse
  recovery problems}, volume 2033 of {\em Lecture Notes in Mathematics}.
\newblock Springer, Heidelberg, 2011.
\newblock Lectures from the 38th Probability Summer School held in Saint-Flour,
  2008, {\'E}cole d'{\'E}t{\'e} de Probabilit{\'e}s de Saint-Flour.
  [Saint-Flour Probability Summer School].

\bibitem{Kol11}
V.~Koltchinskii.
\newblock {\em Oracle inequalities in empirical risk minimization and sparse
  recovery problems: Lecture notes}.
\newblock Springer, 2011.

\bibitem{LeMe16}
G.~Lecu\'e and S.~Mendelson.
\newblock Learning subgaussian classes: Upper and minimax bounds.
\newblock In S.~Boucheron and N.~Vayatis, editors, {\em Topics in Learning
  Theory}. Societe Mathematique de France, 2016.

\bibitem{LeMe16a}
G.~Lecu\'e and S.~Mendelson.
\newblock Sparse recovery under weak moment assumptions.
\newblock {\em Journal of the European Mathematical Society}, to appear.

\bibitem{LeTa91}
M.~Ledoux and M.~Talagrand.
\newblock {\em Probability in Banach Space}.
\newblock Springer-Verlag, New York, 1991.

\bibitem{LeOl12}
M.~Lerasle and R.I. Oliveira.
\newblock Robust empirical mean estimators.
\newblock {\em manuscript}, 2012.

\bibitem{Mas06}
P.~Massart.
\newblock {\em Concentration inequalities and model selection}.
\newblock Ecole d'\'et\'e de Probabilit\'es de Saint-Flour 2003. Lecture Notes
  in Mathematics. Springer, 2006.

\bibitem{Men15}
S.~Mendelson.
\newblock Learning without concentration.
\newblock {\em Journal of the ACM}, 62:21, 2015.

\bibitem{Men16}
S.~Mendelson.
\newblock Local vs. global parameters – breaking the {G}aussian complexity
  barrier.
\newblock {\em manuscript}, 2015.

\bibitem{Men16a}
S.~Mendelson.
\newblock On aggregation for heavy-tailed classes.
\newblock {\em Probability Theory and Related Fields}, to appear.

\bibitem{Men16c}
S.~Mendelson.
\newblock On multiplier processes under weak moment assumptions.
\newblock {\em Geometric Aspects of Functional Analysis - GAFA Seminar notes},
  to appear.

\bibitem{Men16b}
S.~Mendelson.
\newblock Upper bounds on product and multiplier empirical processes.
\newblock {\em Stochastic Processes and their Applications}, to appear.

\bibitem{Min15}
S.~Minsker.
\newblock Geometric median and robust estimation in {B}anach spaces.
\newblock {\em Bernoulli}, 21:2308–2335, 2015.

\bibitem{NeYu83}
A.S. Nemirovsky and D.B. Yudin.
\newblock Problem complexity and method efficiency in optimization.
\newblock 1983.

\bibitem{Sud69}
V.N. Sudakov.
\newblock Gaussian measures, {C}auchy measures and $\epsilon$-entropy.
\newblock In {\em Soviet Math. Dokl}, volume~10, pages 310--313, 1969.

\bibitem{Tal14}
M.~Talagrand.
\newblock {\em Upper and lower bounds for stochastic processes}, volume~60 of
  {\em Ergebnisse der Mathematik und ihrer Grenzgebiete. 3. Folge. A Series of
  Modern Surveys in Mathematics [Results in Mathematics and Related Areas. 3rd
  Series. A Series of Modern Surveys in Mathematics]}.
\newblock Springer, Heidelberg, 2014.
\newblock Modern methods and classical problems.

\bibitem{TS09}
A.~B. Tsybakov.
\newblock {\em Introduction to nonparametric estimation}.
\newblock Springer Series in Statistics. Springer, New York, 2009.
\newblock Revised and extended from the 2004 French original, Translated by
  Vladimir Zaiats.

\bibitem{vdG00}
S.~{van de Geer}.
\newblock {\em Applications of empirical process theory}, volume~6 of {\em
  Cambridge Series in Statistical and Probabilistic Mathematics}.
\newblock Cambridge University Press, Cambridge, 2000.

\bibitem{vdG16}
S.~{van de Geer}.
\newblock {\em Estimation and Testing Under Sparsity}.
\newblock Springer International Publishing, 2016.

\bibitem{vaWe96}
A.W. {van der Vaart} and J.A. Wellner.
\newblock {\em Weak convergence and empirical processes}.
\newblock Springer-Verlag, New York, 1996.

\bibitem{VaCh74a}
V.N. Vapnik and A.Ya. Chervonenkis.
\newblock {\em Theory of Pattern Recognition}.
\newblock Nauka, Moscow, 1974.
\newblock (in Russian); German translation: {\em Theorie der Zeichenerkennung},
  Akademie Verlag, Berlin, 1979.

\end{thebibliography}

\newpage
\appendix
\section{The distance oracle -- outline of the proof}
Here, we sketch the proof of Proposition \ref{thm:distance-functional}. The proof is almost identical to that of Theorem 3.3 from \cite{Men16a}, and follows the same path as the study of the quadratic component in the proof of Theorem \ref{thm:main-general}.

The first step in the proof involves generating a small-ball estimate
that holds with sufficiently high probability, say $9/10$. If a random
variable $Z$ satisfies $\|Z\|_{L_q} \leq L \|Z\|_{L_2}$ for some
$q>2$, then $Z$ automatically satisfies a small-ball estimate, but
with constants $\rho_0$ and $\kappa_0$ that depend on $L$ and $q$. In
particular, it need not be true that $\rho_0$ is close to $1$, for
example, that $\PROB(|Z| \geq \kappa_0\|Z\|_{L_2}) \geq 0.9$. However,
a combination of the norm equivalence and a Berry-Esseen type argument
suffices to ensure that an average of a small number of independent
copies of $Z$ satisfies such a bound. More accurately, if $\ell$ is an
integer that depends only on $q$ and $L$ and $Z_1,\ldots,Z_\ell$ are
$\ell$ independent copies of $Z$, then
\begin{equation} \label{eq:app-small-ball}
\PROB \left(\frac{1}{\ell}\sum_{i=1}^\ell |Z_i|  \geq \kappa_0\|Z\|_{L_2}\right) \geq 0.9~.
\end{equation}
Moreover, combining \eqref{eq:app-small-ball} with a straightforward application of Chebyshev's inequality, one has that with probability at least $0.8$,
$$
\kappa_0 \|Z\|_{L_2} \leq \frac{1}{\ell}\sum_{i=1}^\ell |Z_i| \leq \kappa_1 \|Z\|_{L_2}
$$
for constants $\kappa_0$ and $\kappa_1$ that depend only on $q$ and $L$.

Now consider a partition of $\{1,...,N\}$ to $k$ blocks $I_1^\prime,...,I_k^\prime$, where each block is of cardinality $\ell$. Let ${\cal M}_{Z,j}=\frac{1}{\ell}\sum_{i \in I_j^\prime} |Z_i|$. It follows from a standard binomial estimate that, with probability at least $1-2\exp(-c_1k)$, at least $0.7 k$ of the random variables ${\cal M}_{Z,j}$ satisfy
$$
\kappa_0 \|Z\|_{L_2} \leq {\cal M}_{Z,j} \leq \kappa_1 \|Z\|_{L_2}~.
$$
Let us apply this observation to our setup: fix $f^* \in \F$ and consider the set
$$
{\rm star}(\F-f^*,0) \cap r S \subset {\rm star}(\F-f^*,0) \cap r D = \F_{f^*,r}
$$
for any fixed $r > d^*$. By the choice of $d^*$, ${\rm star}(\F-f^*,0) \cap r S$
contains an $\eta r$-net ${\cal V}_r$ of cardinality at most
$\exp(c_1k/2)$. Therefore, with probability at least
$1-2\exp(-c_1k/2)$, for every $v \in {\cal V}_r$ there is $J_v \subset
\{1,...,k\}$ of cardinality at least $0.7k$, such that for every $j \in J_v$,
\begin{equation} \label{eq:app-two-sided}
\kappa_0 r \leq {\cal M}_{v,j} \leq \kappa_1 r~.
\end{equation}

Next, for each $v \in {\rm star}(\F-f^*,0) \cap r S$ let $\pi v \in {\cal V}_r$ satisfy $\|v-\pi v\|_{L_2} \leq \eta r$. Suppose that one can show that on a high probability event, for every such $v$ there are at most $k/10$ blocks $I_j^\prime$ on which
$$
{\cal M}_{v-\pi v,j} \geq \kappa_0r/2~.
$$
Then, on that event, and on at least $0.6k$ blocks,
$$
(\kappa_0/2) r \leq {\cal M}_{v,j} \leq 2\kappa_1 r~,
$$
and in particular, the same holds for the median of means ${\rm Med}_\ell(v)$.

This is precisely the isomorphic estimate we require, so far, for elements of ${\rm star}(\F-f^*,0)$ whose $L_2(\mu)$ norm is $r$. Obtaining the isomorphic estimate for elements with a larger $L_2(\mu)$ norm follows because ${\rm star}(\F-f^*,0)$ is star-shaped around $0$ and the required isomorphic estimate is positive homogeneous.

Therefore, to complete the proof of the first part of Proposition \ref{thm:distance-functional} it suffices to show that
\begin{equation} \label{eq:app-indicators}
\sup_{v } |\{j : {\cal M}_{v-\pi v,j} \geq \kappa_0r/2\}| \leq k/10,
\end{equation}
where the supremum is taken in ${\rm star}(\F-f^*,0) \cap r S$. The proof of \eqref{eq:app-indicators} is based on an identical argument to the one we used earlier, in the study of the quadratic component, and the fact that $r > d^*$.

The second claim in Proposition \ref{thm:distance-functional} is that on a high probability event, a one-sided (upper) estimate on ${\rm Med}_\ell(v)$ should hold for any $v \in {\rm star}(\F-f^*,0) \cap r D$. Its proof is almost the same as the one we have just described.

For each $v \in {\cal V}_r^\prime$, with probability at least $0.8$,
\begin{equation} \label{eq:app-one-sided-M}
{\cal M}_{v,j} \leq \kappa_1 \|v\|_{L_2} \leq \kappa_1 r~,
\end{equation}
and thus, with probability at least $1-2\exp(-c_1k)$, at least $0.7 k$ of the random variables ${\cal M}_{v,j}$ are smaller than $\kappa_1 r$. Next, we consider an $\eta r$-net ${\cal V}^\prime_{r} \subset \F_{f^*,r}$, which, by the choice of $r$ is of cardinality at most $\exp(c_1k/2)$.
It follows that, with probability at least $1-2\exp(-c_1k/2)$, \eqref{eq:app-one-sided-M} holds for every $v \in {\cal V}_r^\prime$. The oscillation term is then controlled as we outlined above.

Finally, because $k=c(q,L)N$, the estimates hold with the claimed probability of $1-2\exp(-cN)$ for a constant $c$ that depends only on $q$ and $L$.
\endproof

\end{document}